\documentclass[opre,nonblindrev]{informs_preprint}    

\OneAndAHalfSpacedXI 

\usepackage{natbib}
 \NatBibNumeric
 \def\bibfont{\small}%
 \def\bibsep{\smallskipamount}%
 \def\bibhang{24pt}%
 \bibpunct[, ]{[}{]}{,}{n}{}{,}%

\usepackage[colorlinks=true,breaklinks=true,bookmarks=true,urlcolor=blue,
     citecolor=blue,linkcolor=blue,bookmarksopen=false,draft=false]{hyperref}

\def\EMAIL#1{\href{mailto:#1}{#1}}

\TheoremsNumberedThrough     

\EquationsNumberedThrough    

\usepackage{amsmath}
\usepackage{amssymb}
\usepackage{graphicx}
\usepackage{url}
\usepackage{epsfig}\usepackage{psfrag}
\usepackage{MnSymbol}
\usepackage{subfigure}\usepackage{epstopdf}\usepackage{algorithmic}\usepackage{enumerate}
\usepackage{pdfsync}
\usepackage{color}
\usepackage{subfigure}

\long\def\symbolfootnote[#1]#2{\begingroup\def\thefootnote{\fnsymbol{footnote}}
\footnote[#1]{#2}\endgroup}

\global\long\def\bydef{\stackrel{\triangle}{=}}
\global\long\def\buz{\mathbf{Z}}
\global\long\def\zp{\mathbb{Z}_+}
\global\long\def\bls{\mathbf{s}}
\global\long\def\blw{\mathbf{w}}
\global\long\def\blp{\mathbf{p}}
\global\long\def\blq{\mathbf{q}}
\global\long\def\blv{\mathbf{v}}

\global\long\def\blx{\mathbf{x}}
\global\long\def\bly{\mathbf{y}}
\global\long\def\blu{\mathbf{u}}
\global\long\def\calr{\mathcal{R}}
\global\long\def\R{\mathbb{R}}
\global\long\def\pb{\mathbb{P}}

\global\long\def\ity{\infty}

\global\long\def\nor#1{\left\Vert #1\right\Vert }
\global\long\def\norp#1{\left\Vert #1\right\Vert_{\mathbf{p}}}
\global\long\def\norq#1{\left\Vert #1\right\Vert_{\mathbf{q}}}
\global\long\def\E#1{\mathbb{E}\left(#1\right)}

\global\long\def\rkp{\mathbb{R}_{+}^K}
\global\long\def\simp{\Delta}
\global\long\def\bus{\mathbf{S}}
\global\long\def\proj#1{\sigma_{\simp}\left(#1\right)}

\global\long\def\phif#1{\phi^{F}\left(#1\right)}
\global\long\def\phifa{\phi^{F}}
\global\long\def\phisa{\phi^{S}}
\global\long\def\phir#1{\phi\left(#1\right)}
\global\long\def\phis#1{\phi^{S}\left(#1\right)}
\global\long\def\ceils#1{\left\lceil #1\right\rceil }

\global\long\def\rcol#1{\mathcal{R}_{\left(\cdot,#1\right)}}

\global\long\def\highpt#1{h_{\beta}\left(#1\right)}
\global\long\def\lowpt#1{l_{\beta}\left(#1\right)}

\global\long\def\blsw{\blw^{*}}
\global\long\def\bltw{\tilde{\blw}}
\global\long\def\blts{\tilde{\bls}}
\global\long\def\gamap#1{\Gamma\left(#1\right)}

\global\long\def\calc{\mathcal{C}}
\global\long\def\rreal#1{\overline{R}\left(#1\right)}

\begin{document}

\RUNAUTHOR{Mannor and Xu}

\RUNTITLE{Rate-Optimal Control for Resource-Constrained Branching Processes}

\TITLE{Go Viral, or Not: Rate-Optimal Control for Resource-Constrained Branching Processes}

\ARTICLEAUTHORS{%
\AUTHOR{Shie Mannor}
\AFF{Department of Electrical Engineering, Technion, Israel, 32000\\ \EMAIL{shie@ee.technion.ac.il}}
\AUTHOR{Kuang Xu}
\AFF{Laboratory for Information and Decision Systems, MIT, Cambridge, USA, 02139\\ \EMAIL{kuangxu@mit.edu}}
} 

\ABSTRACT{%
We propose and analyze a new class of controlled multi-type branching processes with a per-step linear resource constraint, motivated by potential applications in viral marketing and cancer treatment. We show that the optimal exponential growth rate of the population can be achieved by maintaining a fixed proportion among the species, for both deterministic and stochastic branching processes. In the special case of a two-type population and with a symmetric reward structure, the optimal proportion is obtained in closed-form. 

In addition to revealing structural properties of controlled branching processes, our results are intended to provide the practitioners with an easy-to-interpret benchmark for best practices, if not exact policies. As a proof of concept, the methodology is applied to the linkage structure of the 2004 US Presidential Election blogosphere \cite{AG05}, where the optimal growth rate demonstrates sizable gains over a uniform selection strategy, and to a two-compartment cell-cycle kinetics model for cancer growth, with realistic parameters \cite{DH06}, where the robust estimate for minimal treatment intensity under a worst-case growth rate is noticeably more conservative compared to that obtained using more optimistic assumptions. 
}




\maketitle


\section{Introduction}
\label{sec:intro}
We study \emph{growth-rate maximizing} control strategies, and the associated maximum growth rates, for a new class of multi-type branching processes, where the size of the reproductive population in each round must satisfy a \emph{linear inequality} with respect to the current population. The mathematical model is inspired by, but not limited to, potential applications in marketing, where a firm is interested in maximizing the growth of its active customer base subject to a budgetary constraint, and in cancer treatment, where a doctor would like to estimate the minimal dosage level required to guarantee that the number of the cancer cells diminishes. We begin by describing two motivating examples. For simplicity, all dynamics therein are assumed to be deterministic for the moment, and the introduction of mathematical formalisms is postponed to Section \ref{sec:model}.

%

\subsection{Example 1: Viral Marketing under Budgetary Constraints} 
\label{sec:market}
Consider a marketing campaign where a firm is trying to promote sales among a heterogeneous population over multiple time slots, $t\in \{1,2,\ldots\}$. The underlying population of potential customers is assumed to be of $K$ types, distinguished by their social connectedness and spending patterns.  All customers within the same type are assumed to be identical. The set of \emph{active customers} at slot $t$ is represented by a $K$-dimensional vector $\blw(t)=(\blw_{1}(t),\cdots,\blw_{K}(t))^{\top}$, where $\blw_i(t)$ is the number of type-$i$ customers. During each time slot, a type-$i$ customer generates a revenue of $\blp_i$. Hence, the total revenue for the firm during slot $t$ is given by the weighted $L_1$ norm: $\nor{\blw(t)}_\blp$, where $\nor{\blw}_\blp \bydef \sum_{i=1}^K \blp_i\cdot\blw_i$. At the end of each time slot, the firm strives to expand its customer base by offering  \emph{promotional coupons} to the active customers. A coupon sent to a customer of type $i$ incurs a cost of $\blq_i$.\footnote{We will assume that each customer receives at most one coupon.} If offered a coupon, a customer of type $j$ is expected to recruit $\calr_{(i,j)}$ new customers of type $i$ in the next time slot. A customer who is not offered a coupon does not recruit anyone. We assume that all customers remain active only for one time slot.\footnote{The model can also capture the situation where a customer remains active as long as she receives a coupon. This can be done by adding $1$ to the value of $\calr_{(i,i)}$ for all $1\leq i\leq K$.} Note that the $K\times K$ matrix $\calr$, whose entry on the $i$th row and $j$th column is $\calr_{(i,j)}$, captures the (social) structure of how the underlying populations interact. Denoting by $\bls_i(t)$ the number of type-$i$ active customers who receive coupons at slot $t$, the vector of active customers in $t+1$ is given by 
\begin{equation}
	\blw(t+1) = \calr \bls(t), \nonumber
\end{equation}
Due to budgetary constraints, however, the firm may only afford to send coupons to a \emph{subset} of the active customers. Assuming the firm invests a fraction $\beta$ of its total revenue into marketing in each slot, where $\beta\in[0,1]$, 
\begin{equation}
\label{eq:linConstr}
	\norq{\bls(t)} \leq \beta \norp{\blw(t)},
\end{equation}
and defining the \emph{growth rate} of $\blw(t)$ as\footnote{Alternatively, one could use any weighted $L_1$ norm here (e.g., weighted by the revenue per-individual), since the exponential rate of growth of the $L_1$ norm is insensitive to changes in the weights.}
\begin{equation}
	\alpha=\limsup_{t\rightarrow\ity}\frac{1}{t}\ln\nor{\blw(t)},
	\label{eq:grate}
\end{equation}
the firm may be interested in the following questions.
\begin{enumerate}
	\item What is the \emph{minimum value of $\beta$} required to guarantee a \emph{positive} growth rate of $\nor{\blw(t)}$?
	\item Given a fixed $\beta$, \emph{how many customers in each type} should the firm target (i.e., send coupons to) in each time slot, in order to maximize the growth rate of $\blw(t)$, which is the same as that of the total revenue?
\end{enumerate}

\subsection{Example 2: Robust Intensity Benchmark for Cancer with Cellular Heterogeneity} 
\label{sec:tumor}

In a cancer treatment, the goal may be to minimize, instead of maximize, the growth of the tumor cell population. It is well known that the cellular composition of a cancer can exhibit a high degree of heterogeneity, where multiple types of cancer cells co-exist, with distinct genotypes and reproductive dynamics \cite{DH06,He84,PC07,DF11,DDHM12}.\footnote{One such example of cancer heterogeneity called the active-quiescent cell model, is examined in Section \ref{sec:cancer_study}.} Suppose a cancer therapy (e.g., radiation therapy or chemotherapy) is composed of separate rounds of treatment. Denote by $\blw(t)=\left(\blw_1(t),\ldots,\blw_K(t)\right)^\top$ the state of a cancer with cellular heterogeneity at the beginning of round $t$, where $\blw_i(t)$ is the number of cells of type $i$, and by $\bls_i(t)$ the number of type-$i$ cells at the end of the round, after a treatment has been administered. Similar to the previous example, the treatment's impact on the evolution of the cancer cells can be captured by a controlled multi-type branching process,\footnote{A justification for using a discrete-time branching process to model the cancer growth is provided in Appendix \ref{app:cancer}.} such that $\blw(t+1) = \calr \bls(t)$, $t\in \zp$, where $\calr$ is a $K\times K$ matrix, and  $\calr_{(i,j)}$ is the number of type $i$ cells produced by a cell of type $j$ via mitosis \cite{DH06}, metastasis \cite{DDHM12}, or other mechanisms. 

Unfortunately, unlike the marketing example where the composition of targeted customers can be chosen by the firm, the control decisions here may be in the hands of Nature, or, to be more conservative, \emph{an adversary}. This is because the precision provided by conventional measurement procedures may be insufficient in identifying the composition of the cell types \cite{DH06}, and, even with perfect measurement, it may be impossible to target precisely the desirable composition of cells in a treatment. 

Given the uncertainty, can we obtain a \emph{robust} estimate of the minimum level of treatment intensity that will guarantee certain rate of decrease for the cancer? In particular, assume that the measurement is capable of estimating the \emph{gross size} of the cancer, $\nor{\blw(t)}$, and that the treatment intensity can be tuned so that at least a fraction $1-\beta$ of the cancer cells are exterminated in each round, for some $\beta\in (0,1)$, i.e.,
\begin{equation}
	\nor{\bls(t)} \leq \beta \nor{\blw(t)}.
	\label{eq:cancercons}
\end{equation}
We would like to know, for example, what the minimum value of $1-\beta$ should be so that the growth rate of $\nor{\blw(t)}$ is negative, when $\bls(t)$ is chosen \emph{arbitrarily} from the ones satisfying Eq.~\eqref{eq:cancercons}.  This would tell us the required treatment intensity under a \emph{worst-case} assumption for the pattern of cancer growth. 

\subsection{Remarks on Modeling Assumptions and Objectives} In both the marketing and cancer treatment examples, we have simplified the problem by assuming that the parameter $\beta$ is fixed for all $t$, and that the quantity of interest is the asymptotic growth rate, $\alpha$. One may argue that the problems should be solved over a finite time horizon, and the parameter $\beta$ should be allowed to vary in order to achieve better performance. To this end, the results in this paper are best treated as \emph{first-order} guidelines for solving the problem, rather than exact policies. For instance, while a doctor certainly should not commit to keeping a constant fraction of cell extermination over multiple treatment periods, the maximally allowed value of $\beta$ may still provide her with a simple-to-interpret benchmark for making prescriptions. Finally, our use of a linear constraint is natural for many linear cost functions, such as population size, but maybe ill-fitted for the non-linear ones. 

\subsection{Overview of the Paper}
\label{sec:cont}

We study growth-rate maximizing control strategies for a class of multi-type branching processes with linear resource constraints. The name ``branching'' draws an analogy from the biology literature, in which an individual, if ``chosen'' to reproduce, gives birth to a collection of offspring. While control policies for (multi-type) branching processes have been considered in the literature, a distinguishing feature of our work is that the sub-population $\bls(t)$ chosen to reproduce must obey a \emph{linear inequality} with respect to the current population $\blw(t)$, in the form of $\norq{\bls(t)}\leq \norp{\blw(t)}$.

The majority of our analysis will be focused on a relaxed version of the original integer-valued process, where we allow the number of individuals in each population to take on \emph{continuous values}. As we shall see in the sequel, the use of a continuous population profile allows us to leverage the topological properties of the Euclidean space $\R^K$ and derive results using techniques from the theory of average-cost Markov decision processes (MDP). In particular, our main result shows that the optimal growth rate can always be achieved by maintaining an \emph{optimal mixture} among different population types (Theorem \ref{thm:detcont}). While it may appear that we have made the model more unrealistic by allowing for continuous population profiles, we will show that the optimal growth-rate achievable for an integer-valued {\bf stochastic branching process} (where the size of offspring of each individual is random) coincides with its deterministic continuous counterpart. As a result, the optimal growth rate and its corresponding control strategy carry over to the more realistic stochastic case, where the population is integer-valued, and the reproduction process is random. 

The remainder of the paper is organized as follows. A mathematical formulation of the problem as an average-cost Markov decision process is given in Section \ref{sec:model}. We summarize our main results, along with their interpretations, in Section \ref{sec:result}. Section \ref{sec:DATA} contains two numerical case studies. The first one is based on the social connection data from the 2004 US Presidential Election blogosphere \cite{AG05}, where we show that the optimal strategy can provide a sizable increase in growth rate against a more naive version which selects uniformly across all population types. The second study examines a model of cancer heterogeneity resulted from cell-cycle kinetics, where the cancer cells are divided into an \emph{active} and a \emph{quiescent} compartments. Using a growth model and parameter values in \cite{DH06}, we show that the robust estimate for minimal treatment intensity can be lower than that of a more optimistic estimate, which assumes a uniform extermination rate across all types of cells. The remaining sections are devoted to the proofs, and some concluding remarks are given in Section \ref{sec:concl}.

\subsection{Related Work}
%

Classical (uncontrolled) multi-type branching processes have been extensively studied in the field of applied probability \cite{Har89}, where the special case of a single-type population is known as the Galton-Watson process. The idea of applying control to influence the growth of the branching process dates back to the work of Sevastyanov and Zubkov \cite{SZ74}, and has become an area known as controlled branching processes or Markov population dynamics \cite{Pli77, RW82,GP02}. In this literature, the controller is typically restricted to choosing from a set of actions in each round which would then influence the reproductive behavior of the individuals involved, and one is interested in deriving optimal control laws that maximize the rate of growth of the population (or other notions of rewards) over time. For instance, reference \cite{RW82} shows the existence of a rate-optimal stationary policy in a setting where a control action can be selected independently for each individual. There are two main differences between our approach and this line of work. First, we do not assume control actions for each type (or individual) can be chosen independently from one another, which we believe to be difficult to implement for certain applications especially when the size of the underlying population is large, and second, the \emph{linear resource constraint} imposed on the reproductive population in our model; as will become clear in the sequel, the addition of resource constraints brings analytical challenges that require a different set of tools.

Continuous-time multi-type branching processes have been applied to modeling the growth of cancer with heterogeneity, both in stochastic \cite{DF11,DDHM12} and deterministic \cite{DH06} settings. Most results in this area focus on understanding the dynamics of cancer growth without control (e.g., \cite{DF11,DDHM12}), or under certain specific treatment policies (e.g., \cite{DH06}). The authors of \cite{SPK96} use the theory of optimal control to characterize continuous-time optimal policies that maximize chemotherapy dosage subject to a constraint on a minimum level of bone marrow cell re-generation, where the growth dynamics of the bone marrow cells is captured using the same two-compartment cell-cycle kinetics model used in \cite{DH06}, as well as in the examples of the current paper. In addition to differences in models and objectives, the nature of results in \cite{SPK96} differs from ours in the following two aspects: $1)$ the continuous-time policies do not appear to extend easily to treatments administered over multiple discrete rounds and \emph{vice versa}, and $2)$ the characterization of optimal policies in \cite{SPK96} is given in terms of a solution to a set of differential equations, which provide limited analytical insights or tractability, but may be more relevant for therapies over a small time horizon than the metric of asymptotic growth rate. 

In other application domains, (uncontrolled) branching processes have been used to study the spread of marketing messages  \cite{LBEW10}, information \cite{IM09}, and diseases \cite{Vaz2006} in a (social) network, while there has been little systematic understanding of the impacts of control strategies in such settings. The idea of targeting a subset of the current audience to maximize the spread of information on a network has been studied in the algorithms community \cite{RD02,KKT03}, but under rather different models and dynamics, where the controller is assumed to know the topology of the underlying graph, and the results are typically focused on developing approximation algorithms for NP-hard optimization problems \cite{KKT03}. 

The methodologies used in this paper fall under the umbrella of the theory of Markov decision processes. The reader is referred to \cite{Mai68, Ross83}, and the references therein, for a general introduction to the subject.

\section{Problem Formulation}
\label{sec:model}

\subsection{Notation} We will borrow terminology from the literature of multi-type branching processes (cf. \cite{Har89}). We will refer to the vector $\blw(t)$, whose $i$th coordinate is the number of individuals of type $i$ in time slot $t$, as the {\bf population profile}, and the vector $\bls(t)$, whose $i$th coordinate is the number of individuals of type $i$ allowed to reproduce at slot $t$, as the {\bf reproductive sub-population}, or sub-population, for short. The matrix $\calr$ is referred to as the {\bf reproduction matrix}, because it encodes all the information about the reproductive capabilities of each type of individuals. 

Denote by $\nor{\blx}_\blq$ the weighted $L_1$ norm of $\blx$, with the weights given by $\blq$, and the ordinary unit-weight $L_1$ norm of $\blx$ is simply $\nor{\blx}$.  Denote by $\simp$ the $K$-dimensional simplex, 
\begin{equation}
\simp=\left\{ \blx\in\rkp:\|\blx\|=1\right\},
\end{equation}
and by $\proj{\blx}$ the \emph{scaling of $\blx$ onto $\simp$}, i.e. 
\begin{equation}
\proj{\blx}=\frac{\blx}{\nor{\blx}}.
\end{equation}
We will refer to, $\proj{\blw}$, the scaled version of a population profile $\blw$, a {\bf population mixture}. 

Finally, for any two vectors $\blx,\bly \in\R_{+}^{K}$, we write $\blx\succeq\bly$ (read: ``$\blx$ dominates $\bly$''), if $\blx_{i}\geq\bly_{i},\,\forall1\leq i\leq K$.

\subsection{Dynamic System: REAL} 
\begin{figure}[h]
\centering 
\includegraphics[scale=0.51]{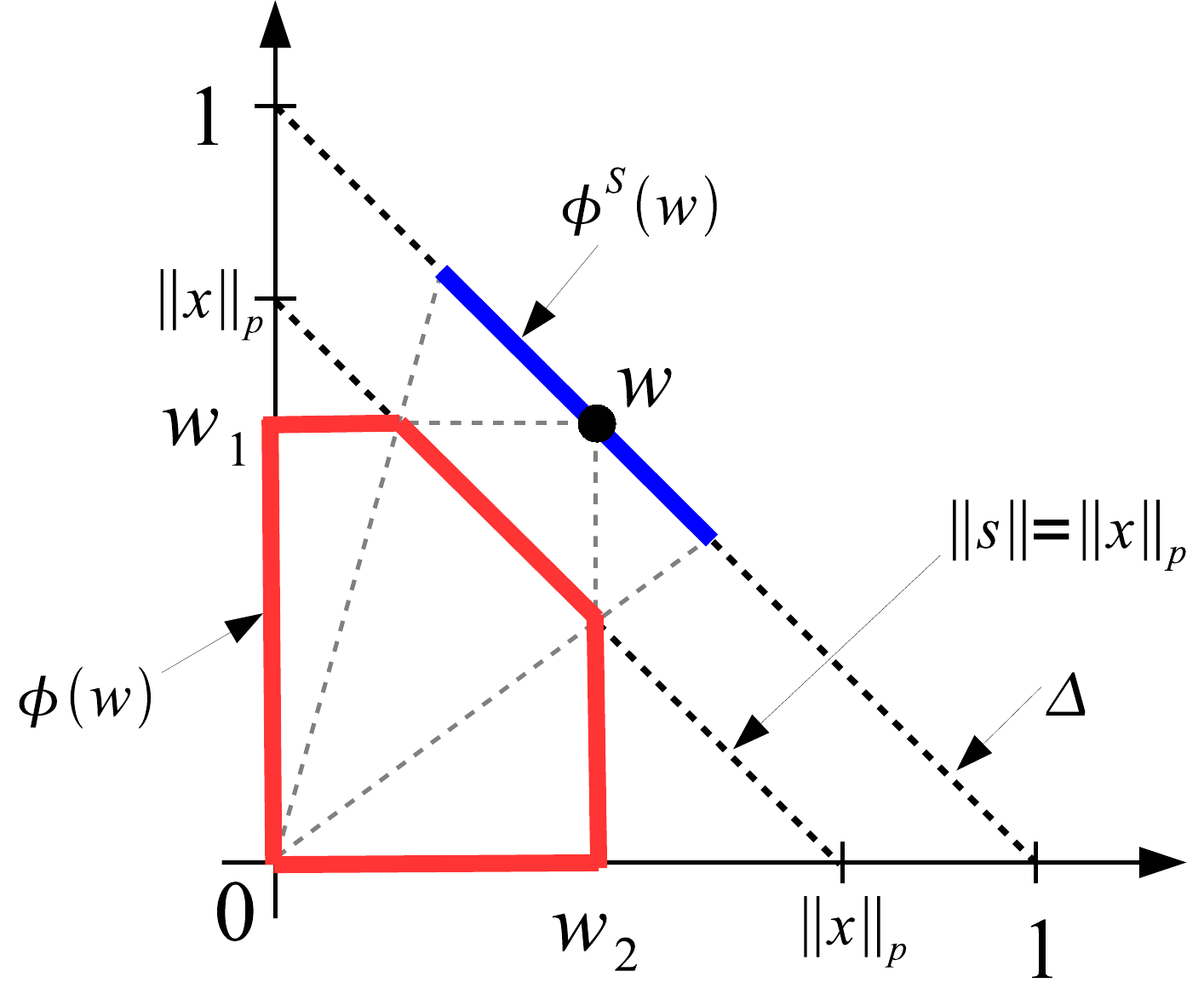} \caption{A geometric view of quantities in REAL and SIM, with $K=2$, where $\blw_i$ is the $i$th coordinate of $\blw$. The feasible action set in REAL, $\phir{\blw}$, is the polyhedron enclosed by the solid red lines (including the interior), and the feasible action set in SIM (defined in Section \ref{sec:REALtoSIM}), $\phis{\blw}$, is marked by the solid blue line segment.}
\label{fig:2dcase1} 
\end{figure}

The following dynamic system formalizes the examples given in Section \ref{sec:intro}, and will serve as the main object of study for this paper.
\begin{definition}
REAL\footnote{REAL is an acronym indicating that the state space of the dynamic system is $\R^K_+$.} is a discrete-time dynamic system given by
\begin{enumerate}
\item \textbf{States}: $\blw(t)\in \rkp$, $ t \in \zp$. 
\item \textbf{Actions}: choose a reproductive sub-population $\bls(t)\in\phir{\blw(t)}$,
where
\begin{equation}
\label{eq:phir}
\phir{\blw}=\left\{ \bls \in\R_{+}^{K}:\nor{\bls}\leq \norp{\blw}\mbox{ and } \bls \preceq  \blw \right\}, 
\end{equation}
Figure \ref{fig:2dcase1} gives an illustration of $\phir{\blw}$ for $K=2$.
\item \textbf{Transition:} $\blw(t+1)=\calr\bls(t)$.
\item \textbf{Reward per-stage}: $\rreal{\blw(t),\bls(t)}=\ln\frac{\nor{\calr\bls(t)}}{\nor{\blw(t)}}$. 
\end{enumerate}
In particular, the parameters for REAL are $1)$ the reproduction matrix $\calr$, and $2)$ the revenue per-individual $\blp$. We assume that all elements of $\calr$ and $\blp$ are strictly positive. 
\end{definition}

\vspace{-18pt}

\begin{definition}
A {\bf feasible policy} is a sequence of mappings $\pi=\left\{\pi_t\right\}_{t\geq 0}$, such that each $\pi_t$ chooses, for all $\blw \in \rkp$, a sub-population $\bls$ in the set $\phir{\blw}$. The set of all feasible policies is denoted by $\Pi$. We say that $\left\{\left(\blw(t),\bls(t)\right)\right\}_{t\geq0}$ is a {\bf feasible sequence} of REAL, if it can be produced by applying some $\pi \in \Pi$. 
\end{definition}

The dynamic system REAL constitutes a (deterministic) Markov decision process. We now comment on why REAL encompasses the examples in Sections \ref{sec:market} and \ref{sec:cancer_study}. First, we note that the change in the linear constraint on sub-population in the definition of $\phir{\blw}$, from the original $\norq{\bls}\leq\beta \norp{\blw}$ in Eq.~\eqref{eq:linConstr} to $\nor{\bls}\leq \norp{\blw}$ in Eq.~\eqref{eq:phir}, is without loss of generality. Clearly, the parameter $\beta$ can be absorbed by changing the $i$th coordinate of $\blp$ to $\beta\blp_i$. To eliminate the dependence on $\blq$, we argue that one can (globally) scale the $i$th coordinate of $\blw$  and $\bls$ by a factor of $\frac{1}{\blq_i}$ (i.e., one unit in the new coordinate is equal to $\frac{1}{\blq_i}$ units in the original), and change the $i$th coordinate of $\blp$ to $\frac{\blp_i}{\blq_i}$. Note that while this re-scaling procedure changes the value of $\nor{\blw(t)}$ by a constant factor, it does \emph{not} change the asymptotic growth-rate of $\nor{\blw(t)}$ as $t\rightarrow \infty$. Second, it is easy to verify that the growth-rate of $\nor{\blw(t)}$ in any feasible sequence $\left\{\left(\blw(t),\bls(t)\right)\right\}_{t\geq0}$ is equal to the corresponding average reward in the limit as $t \to \infty$, by observing that \[\frac{1}{t}\ln\nor{\blw(t)}=\frac{1}{t}\sum_{i=0}^{t-1}\ln\frac{\nor{\blw\left(i+1\right)}}{\|\blw\left(i\right)\|}=\frac{1}{t}\sum_{i=0}^{t}\rreal{\blw(i),\bls(i)}.\] 
Therefore, achieving an optimal growth rate in the original problem is equivalent to having an \emph{average-reward maximizing control policy} in REAL.


\section{Main Results} \label{sec:result}
Our main results are summarized in this section. 

\subsection{Fixed Point Characterization of Rate-Optimal Strategies}
Control policies for REAL can be in general very complex. However,  our first result shows that the \emph{optimal} growth rate can be achieved by keeping the population along a \emph{single} direction, $\blx^* \in \simp$. This significantly reduces the level of complexity faced by the system operator, from an infinite-horizon decision problem to that of finding a single point on the $K$-dimensional simplex. The proof of Theorem \ref{thm:detcont} is given in Section \ref{sec:optpf}. 
\begin{theorem}
\label{thm:detcont} {\bf (Growth-Rate Maximizing Population Mixture)} Fix any positive reproduction matrix $\calr$
and $\blp>0$. There exists an \textbf{optimal population
mixture}, $\blx^{*}\in\simp$, that achieves the optimal growth rate
$\alpha^{*}$. Specifically, 
\begin{equation}
\blw(t)=e^{\alpha^*t} \blx^{*},\quad t\geq0.\label{eq:fixptdef}
\end{equation}
is feasible in REAL, where $\alpha^*$ is the maximum average reward in REAL. 
\end{theorem}

\textit{Remark}: The optimal population mixture given in Theorem \ref{thm:detcont}
can also be interpreted as a \emph{fixed point} of the controlled branching process, in the following sense: there exists $\bls^{*}\in\phi\left(\blx^{*}\right)$
so that $\blx^{*}=\proj{\calr\bls^{*}}$, in which case the optimal growth rate is simply given by $\alpha^{*}=\ln \frac{\nor{\calr\bls^{*}}}{\nor{\blx^{*}}} =\ln  \nor{\calr\bls^{*}}$.\footnote{Since all entries of $\calr$ are assumed to be positive, one can verify that
the optimal growth rate $\alpha^{*}$ guaranteed by Theorem \ref{thm:detcont}
is achievable regardless of the initial condition $\blw(0)$.
 This is true because the optimal population mixture is reachable, up
to a constant factor, starting from any non-zero initial population.
To see this, a policy can let $\blw\left(t\right)=\left(\beta\calr\right)^{t}\cdot\blw\left(0\right),\,\forall1\leq t\leq K-1,$ by letting $\bls\left(t\right)=\beta\blw\left(t\right),\,1\leq t\leq K-1$,
which guarantees that all entries of $\blw(t-1)$ are non-zero. Now,
let $\bls\left(K-1\right)=c\bls^{*},$ for some constant $c>0$,
we have that $\blw\left(K\right)=c\blx^{*}$. In other words, we were able to reach the optimal mixture, up to a factor of $c$, in $K-1$ steps.} 


\subsection{Two-Dimensional Case with Symmetric Revenue}

Theorem \ref{thm:detcont} states that the the optimal growth rate can be achieved via a single population mixture. This existence result, however, does not directly lead to an efficient way of computing the optimal mixture. It turns out that the result can be strengthened if we restrict ourselves to the two-dimensional case ($K=2$), and with a symmetric revenue-per-individual $\blp = \left(\beta,\beta\right)$ for some $\beta\in (0,1)$. Here, the sub-population satisfies
\begin{equation}
\label{eq:2dconstr}
	\nor{\bls(t)} \leq \beta\nor{\blw(t)},
\end{equation}
or, in words, that we are allowed to keep up to a fraction $\beta$ of the total population in each slot. 

In this case, we are able to obtain an explicit characterization of the optimal population mixture in terms of a (unique) fixed point of a \emph{greedy policy}, where the decision maker strictly favors the specifies that can produce more offspring in a single iteration.\footnote{Unfortunately such a characterization does not seem to generalize easily to higher dimension ($K\geq 3$). See discussion preceding the proof of Theorem \ref{thm:k2} in Section \ref{sec:k2pf}.}  The proof of Theorem \ref{thm:k2} is given in Section \ref{sec:k2pf}. 
\begin{theorem}
\label{thm:k2} When $K=2$, the optimal population mixture, $\blx^{*}$, and the associated sub-population,
$\bls^{*}$, are determined by the solution to the following fixed point equations:
\begin{eqnarray}
\label{eq:2dfix1}
\bls_{1}^{*} & = & \min\left\{ \frac{1}{\beta}\blx_{1}^{*},1\right\},\, \bls^*_2 = 1-\bls^*_1,\\ 
\blx^{*} & = & \proj{\calr\bls^{*}}.
\label{eq:2dfix2}
\end{eqnarray}
\end{theorem}
Theorem \ref{thm:k2} provides a means of directly calculating $\blx^{*}$ in the case of $K=2$. In particular, after some elementary algebra, the corollary below follows from solving the fixed point equations \eqref{eq:2dfix1} and \eqref{eq:2dfix2}.
\begin{corollary} \label{cor:k2}
\textbf{(Expressions of Optimal Population Mixture)} Let $\calr=\left(\begin{array}{cc}
a & b\\
c & d
\end{array}\right)$$.$ 
\begin{enumerate}
\item If $\beta<\frac{a}{a-c}$, then $\blx_{1}^{*}=\frac{a}{a-c}$.
\item If $\beta\geq\frac{a}{a-c},$ let $y_{1},y_{2}$ be the two solutions
to the quadratic equation
\[
y^{2}+\frac{a-b\left(\beta+1\right)-\beta d}{\left(b+d\right)-\left(a+c\right)}y+\frac{b}{\left(b+d\right)-\left(a+c\right)}=0.
\]
Then $\blx_{1}^{*}=\max\left\{ y_{i},i=1,2:y_{i}\in\left[0,1\right]\right\}$. 
 
\end{enumerate}
\end{corollary}

\subsection{Stochastic Branching Processes}

We now extend Theorem \ref{thm:detcont} to a more realistic setting, where the population profiles are integer-valued, and the number of offspring produced by a selected individual is \emph{random}. Denote by $\buz(t)=(\buz_{1}(t),\buz_{2}(t),\cdots,\buz_{K}(t))^{\top}\in\zp^{K}$ the population profile at time $t$, where $\buz_{i}(t)$ represents the number of individuals
of type $i$. If a type $i$ individual is selected to reproduce, the set of its off-springs  $\xi_{i}=(\xi_{1,i},\xi_{2,i},\cdots,\xi_{K,i})^{\top}$, is now a \emph{random vector}, drawn from a distribution for which a moment generating function $f_i$ exists. 
The reproductive sub-population that a controller is allowed to select must satisfy the same linear constraint as in the deterministic case, except that now the sub-population must be integer-valued. Specifically, the set of all feasible sub-populations given a current population profile $\buz$ is given by:
\begin{equation}
	\Phi\left(\buz\right)=\left\{ \bus \in\zp^{K}: \nor{\bus} \leq \norp{\buz},\mbox{ and } \bus \preceq \buz \right\}.
	\label{eq:phiz}
\end{equation}
Once a sub-population $\bus(t)$ is chosen from the set $\Phi(\buz(t))$, the distribution of the population in the next time slot is given by 
\[
\buz(t+1)\stackrel{(d)}{=} \sum_{i=1}^{K}\sum_{z=1}^{\bus_i(t)}\xi^{(i,z)}(t),
\]
where the $\xi^{(i,z)}(t)$ are i.i.d (across $z$ and $t$) and distributed according to $\xi_{i}$.

Let $\calr_{(i,j)}=\E{\xi_{i,j}}$ for all $i,j\in\left\{1,\ldots,K\right\}$. The following theorem states that the optimal growth rate associated with the deterministic branching process with reproduction matrix $\calr$ essentially dictates the growth behavior of the stochastic branching process $\buz(t)$. To avoid trivial scenarios of non-extinction, we will assume that
there is a positive probability that an individual will produce zero
offspring, i.e., 
\begin{equation}
\pb\left(\xi_{i,j}=0,\,\forall i \right)>0,\quad\forall1\leq j \leq K.\label{eq:ass1}
\end{equation}

\begin{theorem}
\label{thm:probablistic} Fix any $\blp\succ 0$ and moment generating functions $\left\{f_i\right\}_{1\leq i \leq K}$ for $\left\{\xi_i\right\}_{1\leq i \leq K}$. Let 
\[
\calr_{\left(i,j\right)}=\E{\xi_{i,j}},
\] 
and $\alpha^{*}$ be the optimal growth rate of the corresponding deterministic branching process with reproduction matrix $\calr$, given by Theorem \ref{thm:detcont}. 
\begin{enumerate}
\item If $\alpha^{*}<0$, the stochastic branching process
becomes extinct ($\buz(t)=0$) in finite time with probability one, under any control policy.
\item If $\alpha^{*}>0$, there exists a policy,
$\pi$, under which
\begin{enumerate}
\item Given a sufficiently large initial population, the process explodes
with positive probability, i.e.,
\[
\pb_{\pi}\left(\limsup_{t\rightarrow\infty}\nor{\buz\left(t\right)}=\infty\right)>0.
\]
\item Conditioning on explosion, the growth rate induced by $\pi$ is optimal, i.e.,
\[
\limsup_{t\rightarrow\infty}\frac{1}{t}\ln\nor{\buz\left(t\right)}=\alpha^{*},
\]
almost surely. 
\end{enumerate}
\end{enumerate}
\end{theorem}

\proof{Proof.}
See Section \ref{sec:probpf}. \Halmos\endproof

\section{Numerical Examples}

We present two numerical examples in this section. For a process $\left\{\blw(t)\right\}_{t\in \zp}$ with an associated growth rate $\alpha$, define the \emph{growth factor}, $\kappa$, as
\begin{equation}
	\kappa = e^\alpha,
	\label{eq:kappa_def}
\end{equation}
In other words, $\nor{\blw(t)}=\mathcal{O}\left(\kappa^t\right)$ as $t \to \infty$. We shall use $\kappa$ instead of $\alpha$ in most of our plots for the ease of visualization. 

\subsection{2004 Presidential Election Blogosphere}\label{sec:DATA}
Recall our marketing example in Section \ref{sec:market}. Just as a proof of concept, let us pretend for a moment to be a firm selling political gadgets to a polarized population consisting of only \emph{Liberals} and \emph{Conservatives}. For simplicity, the revenue-per-customer and the cost of each coupon are both $\$1$. Therefore, the only parameter of the model is the fraction of revenue that we are willing to invests in marketing, $\beta$, and the feasible set of sub-population $\bls(t)$ is given by $\nor{\bls(t)}\leq \beta\nor{\blw(t)}$. 

We will study the optimal growth rate of this deterministic branching process using the underlying connectivity structure of the political blogosphere during the period of the 2004 US Presidential Election campaigns \cite{AG05}. In particular, the article documents the linkage structure between the liberal and conservative blogs during the campaign period. Labeling Liberals as type $1$ and Conservatives as type $2$, we will define the reproduction matrix as the following:
\begin{align*}
\calr_{(i,j)} = & \mbox{ average number of links from} \mbox{ a type $j$ blog to a type $i$ blog},
\end{align*}
where the average is computed over the 10 highest-ranked political blogs in the US (\cite{AG05}, Table 1). Based on the data, we have\footnote{The entries are rounded to the nearest integer for simplicity.}
\begin{equation}
\label{eq:R2d}
\calr = \left( \begin{array}{cc}
143 & 17 \\
24 & 137 \end{array} \right).
\end{equation}


%

\begin{figure}[ht]
\vspace{-20pt}
\centering
\subfigure[Growth-factor gain between that of the optimal strategy, $\kappa^*=e^{\alpha^*}$, and that produced by a uniform selection strategy, $\kappa_u=e^{\alpha_u}$.]{
	\includegraphics[scale=.28]{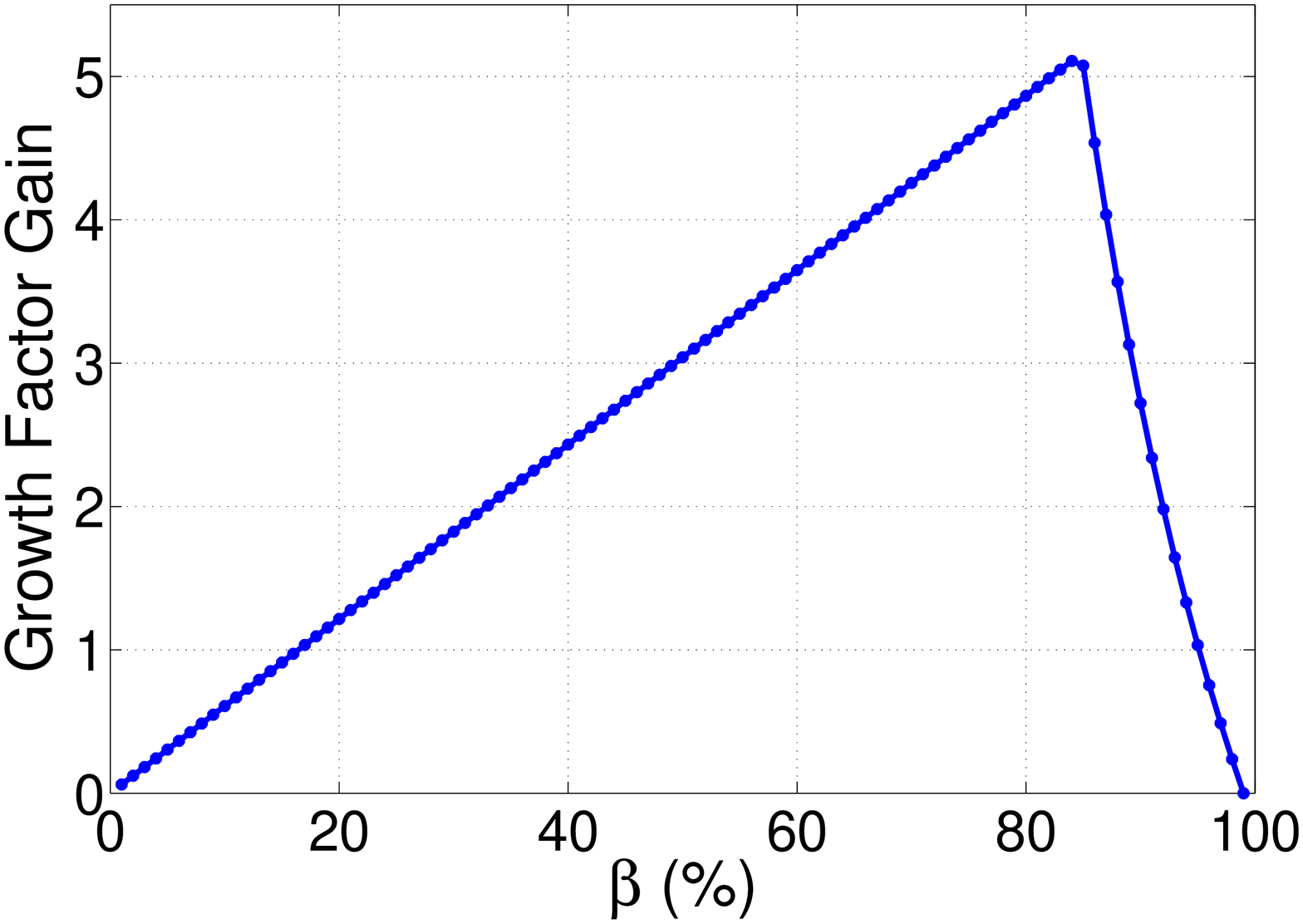}
	\label{fig:1}
}
\subfigure[Fraction of Liberals in the optimal population \newline mixture.]{
	\includegraphics[scale=.28]{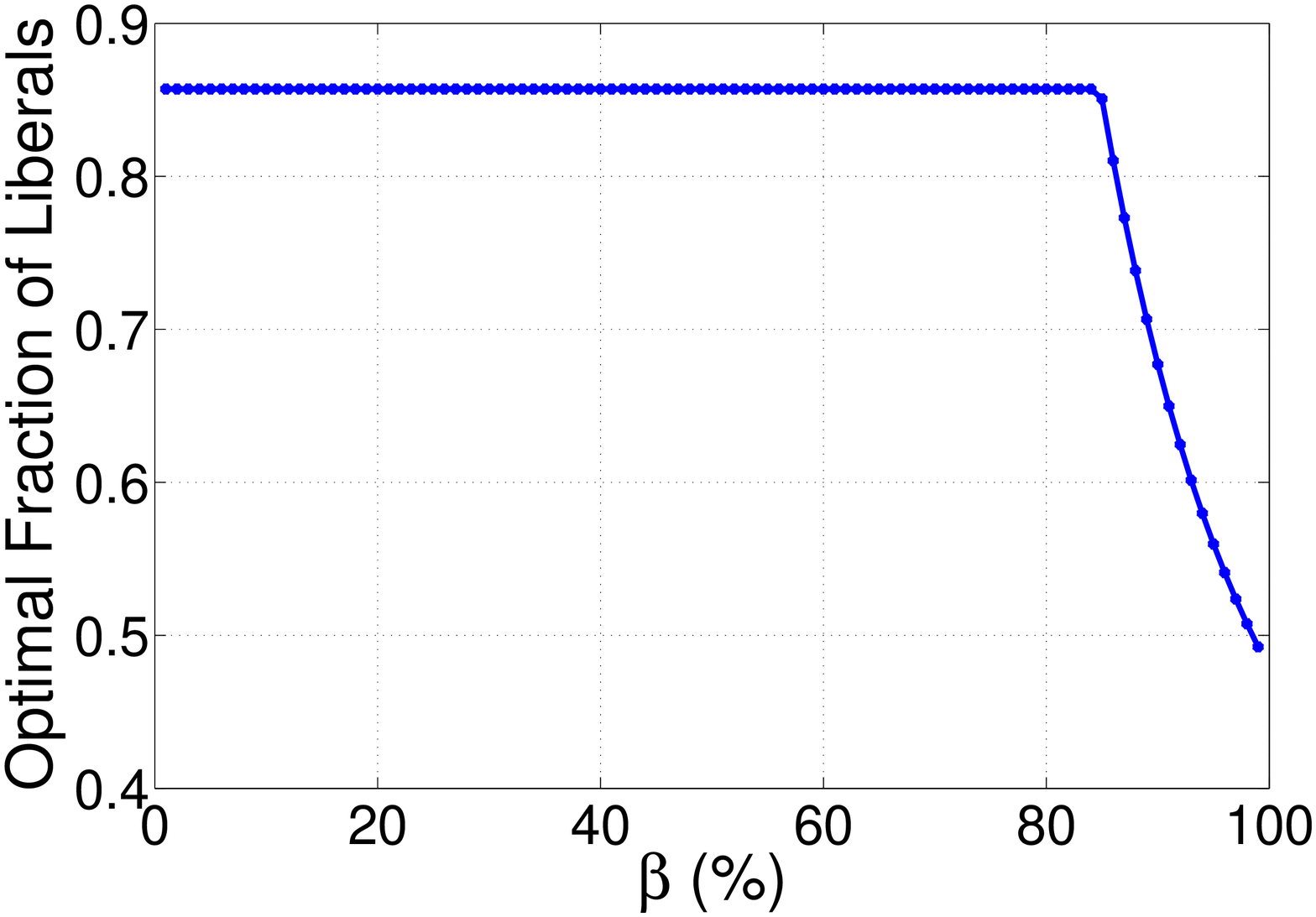}
	\label{fig:2}
}
\label{fig:subfigureExample}
\caption{Example: 2008 Presidential Election Blogosphere.}
\end{figure}

We compare the optimal growth factor (Eq.~\eqref{eq:kappa_def}), $\kappa^*=e^{\alpha^*}$, to that of a naive \emph{uniform policy}, $\kappa_u = e^{\alpha_u}$, which keeps $\beta$ fraction of each type ($\bls(t)=\beta\blw(t)$). The value of $\kappa^*$ is obtained from the optimal mixture given by Corollary \ref{cor:k2}. It is easy to show that the growth rate obtained from the uniform policy is given by $\alpha_u = \beta \rho\left(\calr\right)$, where $\rho\left(\calr\right)$ is the spectral radius (largest eigenvalue) of the matrix $\calr$. Figure \ref{fig:1} shows the difference $\left(\kappa^*-\kappa_u\right)$ as a function of $\beta$. Note that the two values agree near $\beta=0$, where no individual is selected to reproduce, and near $\beta=100\%$, where all of the current population can be selected. In contrast, when there is a moderate amount of flexibility in choosing which type to favor in reproduction ($\beta \approx 80\%$), the optimal growth rate shows sizable gains over the uniform strategy in terms of growth factor.

Figure \ref{fig:2} shows the composition of the {optimal population mixture}, $\blx^*$. Observe from Eq.~\eqref{eq:R2d} that the Liberals have a \emph{higher column sum} compared to the Conservatives (167 v.s 154). We hence know, from Theorem \ref{thm:k2}, that $\blx^*$ must be a fixed point with respect to a greedy selection strategy that maximizes the fraction of Liberals in the sub-population. When resources are limited (small $\beta$), the fixed point corresponds to having \emph{only} Liberals in the reproductive sub-population. As $\beta$ exceeds the level of $80\%$, there becomes enough resources to target also the Conservative population, \emph{in addition} to the Liberals. Finally, as $\beta\rightarrow 1$, since almost the entire population can be selected to reproduce, the gain of using the optimal fixed point over the naive uniform selection strategy diminishes, and the optimal mixture quickly converges to the eigenvector which corresponds to the largest eigenvalue of $\calr$,  $(0.51,\, 0.49)^{\top}$. 

\subsection{Robust Benchmarks for Cancers with Active and Quiescent Cells}
\label{sec:cancer_study}

We examine in this section a particular case of cancer heterogeneity as a result of cell-cycle kinetics \cite{DH06}. Roughly speaking, the cancer cells are divided into two compartments: active cells and quiescent cells. An active cell is capable of reproduction, and after undergoing a mitosis, it splits into two daughter quiescent cells. A quiescent cell does not reproduce, but it may become active after some time period. This cell-cycle dynamics is captured by a set of ordinary differential equations (ODE) in \cite{DH06}. Assuming that the treatment takes place once every of $3$ weeks\footnote{The period of 3 weeks is fairly arbitrary and should vary depending the application. It was chosen as a rough approximation for a typical spacing between two chemotherapies.}, and labeling the active cells as type $1$ and quiescent cells as type $2$, the solutions to the ODEs in \cite{DH06}, as well as the parameter values used, give rise to a discrete-time multi-type branching process, with a reproduction matrix
\begin{equation}
	\calr = \left( \begin{array}{cc} 0.75 & \,0.4674 \\ 1.2864  & \, 0.9258 \end{array} \right),
	\label{eq:cancer_matrix}
\end{equation}
and the reader is referred to Appendix \ref{app:cancer} for a derivation of the discrete-time model from the continuous-time ODEs, as well as the details of how the entries of $\calr$ are computed from the parameters in \cite{DH06}.

\begin{figure}[h]
\vspace{-20pt}
\centering
\includegraphics[scale=.35]{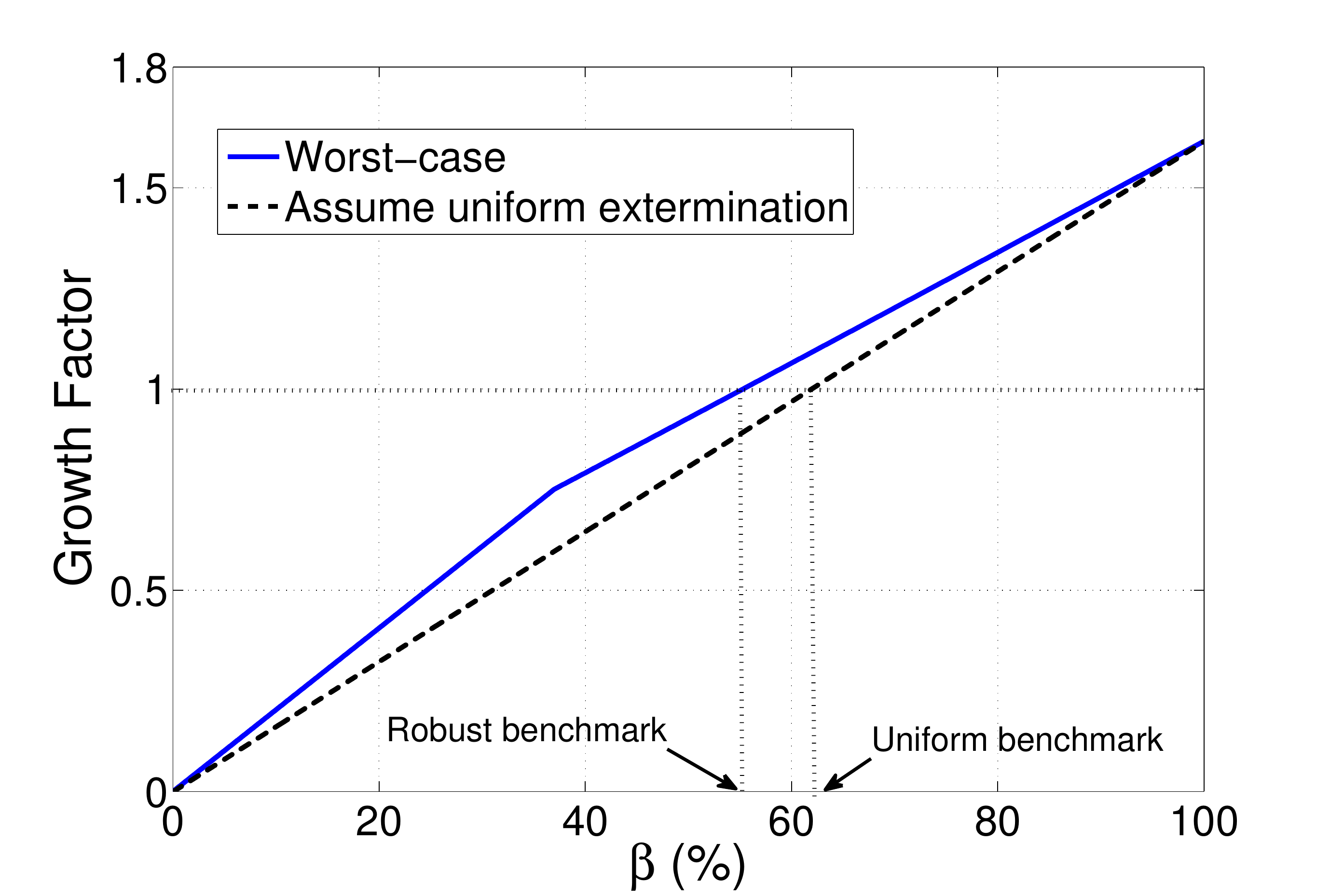}
\caption{Growth rates of the tumor under a uniform extermination assumption versus the worst-case scenario. The time interval is three weeks.}
\label{fig:3}
\vspace{-20pt}
\end{figure}

Figure \ref{fig:3} illustrates the growth factors of $\nor{\blw(t)}$ for $\blp=\blq=(1,1)$.  The worst-case rate is the maximum growth factor, $\kappa^*=e^{\alpha^*}$, induced by the optimal mixture given in Corollary \ref{cor:k2}, where we assume that the composition of the surviving cells after each around of treatment is selected by an adversary. The uniform extermination rate is obtained by assuming that cells of both types are reduced to a fraction of $\beta$ in each round. According to the figure, in order to achieve a diminishing number of cancer cells ($\alpha<0$ or $\kappa<1$), the uniform assumption indicates that at least $37\%$ of the cancer cells should be exterminated per round, whereas the robust estimate is noticeably more conservative, and requires a extermination fraction of $45\%$. 

\section{Proof of Theorem \ref{thm:detcont}}
\label{sec:optpf}
The remainder of the paper is devoted to proofs. Theorem \ref{thm:detcont} is proved in this section, and we begin by outlining the main steps involved. We first define SIM, a dynamic system that evolves on the $K$-dimensional simplex $\simp$, which is induced by scaling the state space of REAL onto $\simp$. Intuitively, SIM tracks the \emph{mixture} associated with the population, as opposed to the actual population vector itself. We then proceed to show that there exists an optimal control strategy for SIM that admits a \emph{fixed point}, which will have implied Theorem \ref{thm:detcont} by the equivalence between SIM and the original system, REAL. To prove the existence of the fixed point, we will use techniques from average-cost Markov decision processes, through the following steps:
\begin{enumerate}
\item We show that the bias function corresponding to the average-reward problem in SIM is \emph{concave} and \emph{continuous} over its domain, $\simp$. 
\item We prove that the induced optimal policy for SIM (a set-valued map)  is continuous and convex-valued.
\item We invoke Kakutani's theorem to establish the existence of a fixed point for the optimal policy. 
\end{enumerate}

Throughout, we will heavily exploit the interplay between the two dynamic systems, REAL and SIM: the compact state space of SIM allows us to leverage stronger topological properties, such as the tightness of a sequence of Lipschitz-continuous functions, while the dynamics of REAL is more intuitive and convenient to work with in deriving bounds on the value and bias functions for the Markov decision processes. 

\subsection{From REAL to SIM}
\label{sec:REALtoSIM}

We show in this subsection that a dynamic system, SIM, whose state space is the $K$-dimensional simplex, suffices in fully describing the dynamics of our original system, REAL. We first state some properties of the action sets of REAL. The proof of Lemma \ref{lem:supadd} is given in Appendix \ref{app:lem:supadd}. 

\begin{lemma}\label{lem:supadd}
\textbf{(Properties of Action Sets)}
Fix any $\blw,\bltw\in\R_{+}^{K}$, and let $\phif{\blw}$ be defined as in Eq.~\eqref{eq:phir}. The following is true.\footnote{See Definition \ref{def:minkowski} in Appendix \ref{app:techprelim} for the definition of the addition of sets.}
\begin{enumerate}
\item {\bf (Scale-invariance)} $a\phir{\blw} = \phir{a \blw}$, for all $a>0$. 
\item {\bf (Superadditivity)}  $\phir{\blw+\bltw}\supset\phir{\blw}+\phir{\bltw}$,
\item {\bf (Convexity)} $\phir{a\blw+(1-a)\bltw} \supset a \phir{\blw}+(1-a) \phir{\bltw}$, for all $a \in [0,1]$.
\end{enumerate}
\end{lemma}

In principle, starting from any population profile $\blw\in\R_{+}^{K}$,
the feasible set of reproductive population is given by $\phir{\blw}$. However, the next lemma
implies that, as far as growth-rate maximizing policies are concerned, it suffices to restrict to a \emph{single facet} on the boundary of $\phir{\blw}$\footnote{This is intuitive, since one should never exclude more individuals from reproducing than necessary, when trying to maximize the growth rate of the process.} , given by\footnote{In Figure \ref{fig:2dcase1},  $\phif{\blw}$ corresponds to the upper-right facet of the polyhedral in red. }
\begin{equation}
	\phif{\blw}\bydef\left\{ \bls\in\R_{+}^{K}:\nor{\bls}=\min\left\{\nor{\blw}, \nor{\blw}_p\right\},\mbox{ and } 0 \preceq \bls \preceq \blw \right\}.
\label{eq:phif}
\end{equation}

\begin{lemma}
\textbf{\label{lem:facet} (Monotonicity of Feasible Sequences)} Fix
any $\blw,\bltw\in\R_{+}^{K}$,
such that $\blw \succeq\bltw$. For any feasible sequence $\left\{\left(\bltw(t),\blts(t)\right)\right\}_{t\geq 0}$ in REAL with $\bltw\left(0\right)=\bltw$, there exists a feasible sequence $\left\{\left(\blw(t),\bls(t)\right)\right\}_{t\geq 0}$ with $\blw\left(0\right)=\blw$ such that $\blw\left(t\right)\succeq\bltw\left(t\right)$ for all $t\geq 0$. 
\end{lemma}

\proof{Proof.} It suffices to show that $\blw \succeq \bltw$ implies that $\phir{\blw}\supset \phir{\bltw}$, from which the claim will follow via a simple induction. To see this, let $\blu=\blw-\bltw$. Since $\blw\succeq\bltw,$ we have that $\blu\in\rkp$, and that $0\in\phir{\blu}$. By the superadditivity property in Lemma \ref{lem:supadd}, we have that $\phir{\bltw}=\phir{\bltw}+0\subset\phir{\bltw}+\phir{\blu}\subset\phir{\blw}$. 
\Halmos\endproof

To see why Lemma \ref{lem:facet} justifies restricting our attention to the facet, choose any $\bls\in\phir{\blw\left(t\right)}$ that does not belong
to the facet $\phif{\blw\left(t\right)}$. It is easy to verify that
there exists some $\hat{\bls}\in\phif{\blw\left(t\right)}$ such
that $\hat{\bls}\succeq\bls$. 
Since all entries of $\calr$ are non-negative, we have that 
$\calr\hat{\bls}\succeq\calr\bls$. 
By Lemma \ref{lem:facet}, this implies that the average reward is
not compromised by choosing the sub-population $\hat{\bls}$ over
$\bls$. By considering only the facets as our action sets, we now formulate a new dynamic system, where the state space is the $K$-dimensional simplex, $\simp$, which yields the same optimal average reward as REAL. 

\begin{definition}
SIM\footnote{SIM is an acronym indicating that the state space of the dynamic system is the $K$-dimensional simplex, $\simp$.} is a discrete-time dynamic system with
\begin{enumerate}
\item \textbf{States}:  $\blw(t)\in \simp$, $\forall t \in \zp$.
\item \textbf{Actions}: choose a reproductive sub-population $\bls(t)\in\phis{\blw}$,
where $\phis{\blw}$ is defined to be the set $\phif{\blw}$ scaled to the simplex (c.f., Figure \ref{fig:2dcase1})
\begin{equation}
\phis{\blw}\bydef\proj{\phif{\blw}}.\label{eq:phis}
\end{equation}
Given the expression for $\phif{\blw}$ (Eq.~\eqref{eq:phif}), we also have that
\begin{equation}
	\phis{\blw} = \left\{ \bls \in\simp: 0\leq\bls_{i}\leq\frac{1}{\min\left\{\norp{\blw}, \nor{\blw}\right\}}\blw_{i},\,\forall1\leq i\leq K\right\}. 
	\label{eq:phisexp}
\end{equation}
\item \textbf{Transition:} $\blw(t+1)=\proj{\calr\bls\left(t\right)}$. 
\item \textbf{Reward per-stage}: $R\left(\blw(t), \bls(t)\right)=\ln \left(\nor{\calr\bls(t)}\right) + \ln\left(  \min\left\{\norp{\blw},1\right\}\right)$.\footnote{The second term, $\ln\left(  \min\left\{\norp{\blw},1\right\}\right)$, accounts for the effect of scaling $\calr\bls\left(t\right)$ onto $\Delta$.}
\end{enumerate}
\end{definition}

We conclude this subsection by stating the following lemma, which implies that REAL and SIM are essentially \emph{equivalent}. It follows directly from Lemma \ref{lem:facet} and the definitions of REAL and SIM. 
\begin{lemma}
\label{lem:SIM-REAL-EQ}
For every feasible sequence $\left\{\left(\blw^S\left(t\right),\bls^S(t)\right)\right\}_{t\geq 0}$ in SIM, there exists a feasible sequence $\left\{\left(\blw^R\left(t\right),\bls^R(t)\right)\right\}_{t\geq 0}$ in REAL, with $\blw^R(0)=\blw^S(0)$, such that 
\[R\left(\blw^R(t),\bls^R(t)\right)=R\left(\blw^S(t),\bls^S(t)\right), \quad \forall t \geq 0.\]
Conversely, for every feasible sequence $\left\{\left(\blw^R\left(t\right),\bls^R(t)\right)\right\}_{t\geq 0}$ in REAL, there sequence $\left\{\left(\blw^S\left(t\right),\bls^S(t)\right)\right\}_{t\geq 0}$, where $\blw^S(t)=\proj{\blw^R(t)}$ for all $t\geq 0$, is feasible in SIM, and satisfies
\[R\left(\blw^S(t),\bls^S(t)\right) = R\left(\blw^R(t),\bls^R(t)\right), \quad \forall t \geq 0.\]
\end{lemma}

\subsection{Markov Decision Processes in SIM}

As was shown in the preceding subsection, the dynamic system SIM is equivalent to the original system, REAL. Therefore, it suffices to focus our attention on the property of the average-reward Markov decision process (MDP) in SIM. The objective of this subsection is to show that the bias function for the average-reward MDP in SIM is \emph{concave} and \emph{continuous}. We will establish this fact in two steps:
\begin{enumerate}
\item  (Propositions \ref{prop:val_convx-1} and \ref{prop:val_lipcont}) Show that there exists $l>0$, such that the value function of the \emph{$\gamma$-discounted-reward} MDP in SIM is {concave} and {$l$-Lipschitz continuous} over $\simp$, for any discount factor $\gamma\in (0,1)$. 
\item (Proposition \ref{prop:adv_val_conv}) Show that the $\gamma$-discounted value function converges uniformly to the average-reward bias function as $\gamma\to 1$, under appropriate normalization. Combined with the previous step, this implies that concavity and continuity properties hold for the average-reward bias function. 
\end{enumerate}

We begin by considering an infinite-horizon discounted problem on SIM. Given a policy $\pi$, the discounted total reward from a state
$\blw$ is defined as $V_{\pi}^{\gamma}\left(\blw\right)=\sum_{t=0}^{\infty}\gamma^{t}R\left(\blw(t),\bls(t)\right)$, where $\gamma\in(0,1)$ is the discount factor, and $\left\{ \left(\blw(t),\bls\left(t\right)\right)\right\} _{t\geq0}$
is the induced sequence by applying policy $\pi$, with
$\blw(0)=\blw$. The \textbf{$\gamma$-discounted value function,}
$V^{\gamma}\left(\cdot\right)$, is defined by 
\begin{equation}
V^{\gamma}\left(\blw\right)=\sup_{\pi\in\Pi}V_{\pi}^{\gamma}\left(\blw\right),\quad\forall\blw\in\simp.\label{eq:vinf}
\end{equation}

The next following result is a form of Bellman's equation in our setting. The proof of is standard, and is given in Appendix \ref{app:lem:discount_Bell} for completeness. 

\begin{lemma}
\textbf{\label{lem:discount_Bell}}
$V^{\gamma}\left(\cdot\right)$ satisfies the fixed point
equation
\begin{equation}
V^{\gamma}\left(\blw\right)=\max_{\bls\in\phis{\blw}}\left[R\left(\blw, \bls\right)+\gamma V^{\gamma}\left(\proj{\calr\bls}\right)\right],\quad\forall\blw\in\simp.\label{eq:disBell}
\end{equation}
In addition, if a stationary policy $\pi^{*}$attains the minimum
in Eq.~\eqref{eq:disBell} for all $\blw\in\simp$, then $\pi^{*}$
is optimal, i.e., $V_{\pi^{*}}\left(\blw\right)=V^{\gamma}\left(\blw\right),\,\forall\blw\in\simp$.
\end{lemma}

\subsubsection{Concavity of Discounted Value Function}

We show in this subsection that the value function for the discounted problem is always concave. The proof for concavity relies a standard coupling argument among multiple alternative state sequences. However, the coupling will be done in the original REAL, as opposed to SIM. This is because in REAL the state sequence alone is sufficient in describing the reward earned at each stage. In contrast, in SIM, due to the scaling onto the simplex in each step, knowing the rewards also requires keeping track of the sequence of actions, which is more complicated to do.  

\begin{proposition}
\label{prop:val_convx-1} {\bf (Concavity of Value Functions)} $V^{\gamma}\left(\cdot\right)$ is a concave function on $\simp$
for all $\gamma\in\left(0,1\right)$. 
\end{proposition}
\proof{Proof.}
Let $\left\{\left(\blw(t),\bls(t)\right)\right\}_{t\geq 0}$ and $\left\{\left(\bltw(t),\blts(t)\right)\right\}_{t\geq 0}$ be two feasible sequences of REAL, with $\blw(0), \bltw(0) \in \simp$. Fixing any $a\in [0,1]$, by the convexity of action sets of REAL (Lemma \ref{lem:supadd}), the sequence $\left\{\left(\blw^a(t),\bls^a(t)\right)\right\}_{t \geq 0}$, defined by
\[ \blw^a(t) = a\blw(t)+(1-a)\bltw(t),\mbox{ and } \bls^a(t) = a\bls(t)+(1-a)\blts(t), \quad \forall t\geq 0, \]
is also feasible. Denote by $R_{\blw}$, $R_{\bltw}$ and $R_a$ the associated $\gamma$-discounted reward in REAL for the sequences $\blw(t)$, $\bltw(t)$ and $\blw^a(t)$, respectively. We claim that
\begin{equation}
\label{eq:concav1}
R_a \geq aR_{\blw} + (1-a)R_{\bltw},
\end{equation}
for all $\gamma\in(0,1)$. By the equivalence of REAL and SIM (Lemma \ref{lem:SIM-REAL-EQ}), this would imply that 
\[V^\gamma\left(a\blw+(1-a)\bltw \right) \geq a V^\gamma_\pi\left(\blw\right) + (1-a) V^\gamma_{\pi'} \left(\bltw\right), \]
for all $a\in [0,1]$, $\blw,\bltw \in \simp$, and $\pi,\pi'\in \Pi$, which in turn yields the concavity of $V^\gamma\left(\cdot \right)$. 
We now show Eq.~\eqref{eq:concav1}. Since $\nor{\blw^a(0)}=\nor{\blw(0)}=\nor{\bltw(0)}=1$, we have
\begin{eqnarray}
\sum_{t=0}^{n-1} \rreal{\blw^a(t),\bls^a(t)} &=& \ln\left(\nor{\blw^a(n)}\right) \geq a\ln\left(\nor{\blw(n)}\right)+(1-a)\ln\left(\nor{\bltw(n)}\right) \nonumber \\
&=& a\left( \sum_{t=0}^{n-1} \rreal{\blw(t),\bls(t)}\right) + (1-a)\left(\sum_{t=0}^{n-1} \rreal{\bltw(t),\blts(t)} \right) \nonumber \\
&=& \sum_{t=0}^{n-1} a\rreal{\blw(t),\bls(t)} + (1-a) \rreal{\bltw(t),\blts(t)}, 
\label{eq:concav2}
\end{eqnarray}
for all $n\geq 1$. Eq.~\eqref{eq:concav1}, and hence our claim, follows by applying the following technical lemma to Eq.~\eqref{eq:concav2}, by letting $a_t =  \rreal{\blw^a(t),\bls^a(t)} $ and $b_t = a\rreal{\blw(t),\bls(t)} + (1-a)\rreal{\bltw(t),\blts(t)}$. 
\begin{lemma}
\label{lem:tech2} Let $\left\{a_t\right\}_{t\geq 0}$ and $\left\{b_t\right\}_{t\geq 0}$ be two sequences such that $\sum_{t=0}^{n-1} a_t \geq \sum_{t=0}^{n-1} b_t,$ for all $n \geq 0$. Then for all $\gamma\in (0,1)$, $\sum_{t=0}^{n-1} \gamma^t a_t \geq \sum_{t=0}^{n-1} \gamma^t b_t$ for all $n\geq 1$. 
\end{lemma}
\proof{Proof.}
See Appendix \ref{app:tech2}.
\Halmos\endproof
This completes the proof of Proposition \ref{prop:val_convx-1}. 
\Halmos\endproof

\subsubsection{Smoothness of Discounted Value Functions}

The concavity of  the value function $V^\gamma(\cdot)$ (Proposition \ref{prop:val_convx-1}) already implies that $V^\gamma(\cdot)$  is continuous in the interior of $\simp$. However, a stronger notion of continuity is needed for our purpose. In this subsection, we will show that $V^\gamma(\cdot)$ is continuous over the entire domain, $\simp$, with a \emph{bounded} Lipschitz coefficient. We begin with the following definition. 

\begin{definition} Fix $l,\epsilon>0$. A function $f:\R^{K}\rightarrow\R$ is $l$-Lipschitz continuous if
$\left|f\left(\blx\right)-f\left(\bly\right)\right|\leq l\nor{\blx-\bly}$, for all $\blx,\bly\in\R^{K}$. A function $f:\R^{K}\rightarrow\R$ is $\epsilon$-locally $l$-Lipschitz continuous if $\left|f\left(\blx\right)-f\left(\bly\right)\right|\leq l\nor{\blx-\bly}$, for all $\nor{\blx-\bly}\leq\epsilon$. 
\end{definition}

The next proposition is the main result of this subsection. 

\begin{proposition}
\label{prop:val_lipcont}{\bf (Smoothness of Value Functions)} There exists $l>0$, such that $V^{\gamma}\left(\cdot\right)$ is $l$-Lipschitz-continuous for all $\gamma\in\left(0,1\right)$.
\end{proposition}

The remainder of this subsection is devoted to the proof of Proposition \ref{prop:val_lipcont}. We first briefly discuss the proof strategy. By definition, establishing the continuity of the discounted value function amounts to showing that, starting from two initial locations $\blw$ and $\bltw$ in $\simp$ that are ``close'' to each other, the optimal discounted rewards from both points are also similar. We will prove this via a coupling argument, that in REAL, the state sequence starting form $\bltw$ can eventually ``emulate'' any sequence starting from $\blw$, by sustaining some penalty in the first few iterations, and hence $V^\gamma(\bltw)$ cannot be too small compared to $V^\gamma(\blw)$. To achieve the emulation, we will exploit a tradeoff intrinsic in the structure of the multi-type branching process: in each step, the policy could choose to keep a \emph{smaller} reproductive population in terms of total size (short-term loss in reward), but in return, it now has more freedom in choosing the \emph{mixture} among the reproductive species (long-term gain). Applying this idea to our context, we will let the sequence starting from $\bltw$ choose a smaller reproductive population than necessary, in exchange for the ability to produce the \emph{same} population mixture as the sequence from $\blw$ in just a few iterations. This tradeoff will be made rigorous in Lemma \ref{lem:cutting}.

Let $\calc$ be the convex hull of all column vectors of $\calr$, scaled to $\simp$,
\[
\calc\bydef conv\left(\left\{ \proj{\rcol 1},\proj{\rcol 2},\ldots,\proj{\rcol K}\right\} \right).
\]
The proof of Proposition \ref{prop:val_lipcont} will be accomplished in two steps:
\begin{enumerate}
	\item (Lemma \ref{lem:lipsinC}) Apply the above-mentioned coupling argument to $\blw,\bltw\in\calc$ that are ``nearby'', and show that the cost of emulation is small. This will prove that $V^\gamma(\cdot)$ is smooth over $\calc$. 
	\item Extend the smoothness of $V^\gamma(\cdot)$ over $\calc$ to the entire domain, $\simp$, thus completing the proof of Proposition \ref{prop:val_lipcont}. 
\end{enumerate}

We begin by defining a subset of the action set in REAL, $\phifa\left(\blw,a\right)$, that captures the notion of ``preserving a smaller reproductive population than necessary'', as was mentioned above. Let
\begin{equation}
\label{eq:phif2}
	\phifa\left(\blw,a\right)  \bydef \left\{\bls \in \rkp: \nor{\bls} = a \min\left\{\norp{\blw},\nor{\blw}\right\} \mbox{ and } 0\preceq \bls \preceq \blw \right\},
\end{equation}
where $\blw\in \rkp$ and $a\in \R_+$, and 
\[\phisa\left(\blw,a\right) \bydef \proj{\phifa\left(\blw,a\right)}.\]
The following consequences are easy to verify 
\begin{enumerate}
	\item $\phifa\left(\blw,a\right)\subset \phir{\blw}$, for all $a\in \left[0,1\right]$, and hence all points in $\phifa\left(\blw,a\right)$ are feasible actions in REAL. 
	\item $\phisa\left(\blw,a \right) \supseteq \phis{\blw}$, for all $a\in \left[0,1\right]$, with $\phisa\left(\blw,1 \right) = \phis{\blw}$. 
\end{enumerate}

The next lemma formalizes the tradeoff between short-term rewards and wider choices of sub-populations mixtures. It states that for
any two points $\blw,\bltw\in\calc$, the set of reachable reproductive sub-population mixtures (scaled to $\simp$) from $\blw$ can also be reached from $\bltw$, if one is willing to keep a smaller sub-population. 
\begin{lemma}
\label{lem:cutting} There exists a constant $\mu>0$, so that
for all $d>0$ and all $\blw,\bltw\in\calc$, $\nor{\blw-\bltw}= d$,
we have
\[
\phis {\bltw}\subset \phisa\left(\blw,(1-\mu d)^2\right).
\]
\end{lemma}

\proof{Proof.} Fix $\bltw \in \calc$. Let $h\left(\blw\right)\bydef \min\left\{\norp{\blw},\nor{\blw}\right\}$. It can be verified from the definition that $\phisa\left(\blw,a\right)$ can be equivalently written as
\[
\phisa\left(\blw,a\right) = \left\{\bls\in \simp: 0 \leq \bls_i \leq \frac{1}{a h\left(\blw\right)} \blw_i, \,\forall 1\leq i\leq K\right\}.
\]
It therefore suffices to show that $\frac{\blw_i}{\left(1-\mu d\right)^2h\left(\blw\right)} \geq \frac{1}{h\left(\bltw\right)}$ for all $1\leq i \leq K$, or, equivalently, that
	\[
	\left(1-\mu d\right)^2 \leq \frac{h\left(\blw\right)}{h\left(\bltw\right)}\cdot \frac{\bltw_i}{\blw_i}, \quad \forall 1 \leq i \leq K.
\]
The following facts are simple consequences of the definitions and assumptions, 
\begin{enumerate}
	\item There exist $a,b >0$ such that $\min_{\blw \in \calc} h\left(\blw\right)>a$, and $\min_{\blw\in \calc, 1\leq i \leq K} \blw_i >b$. 
	\item There exists $c >0$, such that $\left|h\left(\blw\right)-h\left(\bltw\right)\right| \leq c\nor{\blw-\bltw}$, for all $\blw, \bltw \in \R^K$. 
\end{enumerate}
From the above facts, we conclude that there eixsts some $\mu>0$, such that for all $\blw \in \calc$, $\nor{\blw-\bltw}=d$, 
\begin{eqnarray*}
\frac{h\left(\blw\right)}{h\left(\bltw\right)}\cdot \frac{\bltw_i}{\blw_i} &=& \left(1-\frac{h\left(\bltw\right)-h\left(\blw\right)}{h\left(\bltw\right)}\right) \left(1-\frac{\blw_i-\bltw_i}{\blw_i}\right) \geq \left(1-\frac{c}{a} d\right)\left(1-\frac{1}{b}d\right) \\
&\geq& \left(1-\mu d\right)^2,
\end{eqnarray*}
for all $1 \leq i \leq K$, which completes the proof. 
\Halmos\endproof

The next lemma shows that $V^{\gamma}\left(\cdot\right)$ is Lipschitz-continuous
within the set $\calc$. The proof involves a coupling and emulation argument using Lemma \ref{lem:cutting}. 

\begin{lemma}
\label{lem:lipsinC} There exists $l_{1},\epsilon>0$ such that
for all $\blw,\bltw\in\calc,$ $\nor{\blw-\bltw}\leq\epsilon$, we have
that 
\begin{equation}
V^{\gamma}\left(\bltw\right)-V^{\gamma}\left(\blw\right)\leq  l_{1}\nor{\blw-\bltw}.\label{eq:vdisupper1}
\end{equation}
\end{lemma}

\proof{Proof.} Note that by Lemma \ref{lem:localglobal} in Appendix \ref{app:techprelim},
it suffices to show that $V^{\gamma}\left(\cdot\right)$ is $\epsilon$-locally
$l$-Lipschitz continuous. We will show that there exist constants $\mu, \epsilon > 0$, such that for all $\blw,\bltw\in \calc$, $\nor{\blw-\bltw}=d \leq \epsilon$, 
\begin{equation}
	\label{eq:lip0}	\left|V^\gamma\left(\blw\right)-V^\gamma\left(\bltw\right)\right| \leq \left|\ln\left(1-\mu d\right)^3\right|,
\end{equation}
from which the claim would follow by taking the Taylor expansion of $\ln(x)$ around $x=1$. 

Since all entries of $\calr$ are positive, we have $\min_{\blw\in \calc} \norp{\blw} >0$, and hence there exist $\epsilon, \mu_1>0$ such that 
\begin{equation}
\label{eq:lip1}
	\frac{\norp{\bltw}}{\norp{\blw}} \geq 1 - \mu_1 d,
\end{equation}
for all $\nor{\blw-\blv}=d\leq \epsilon$. Fix any $\blw,\bltw \in \calc$, $\nor{\blw-\bltw} = d \leq \epsilon$. By Lemma \ref{lem:cutting}, we have
\begin{equation}
\label{eq:lip2}
	\phis{\bltw} \subset \phisa\left(\blw, \left(1-\mu_2 d\right)^2\right),
\end{equation}
where $\mu_2$ is some positive constant. 

Let $\left\{\left(\bltw\left(t\right),\blts\left(t\right)\right)\right\}_{t\geq 0}$ be a feasible sequence in SIM with $\bltw(0)=\bltw$, and denote by $R_{\bltw}$ the total discounted reward of the sequence: $R_{\bltw} = \sum_{t=0}^\ity \gamma^t R\left(\bltw(t), \blts(t)\right)$. By definition of the value function, we know that $V^\gamma\left(\bltw\right)=R_{\bltw}$ for some feasible sequence starting from $\bltw$. We now employ a coupling argument, by constructing a feasible sequence in REAL (rather than SIM), which starts at $\blw$, that offers a total discounted reward that is ``close'' to $R_{\bltw}$. In particular, consider sequence $\left\{\left(\blw\left(t\right),\bls\left(t\right)\right)\right\}_{t\geq 0}$, generated according to:
\[
\blw\left(0\right)=\blw,
\]
\[
\bls\left(t\right)=\begin{cases}
\left(1-\mu_2 d\right)^2 \blts(t), & \, t=0,\\
\frac{\min\left\{\norp{\blw(t)},\nor{\blw(t)}\right\}}{\nor{\blts(t)}} \blts(t), & t\geq1,
\end{cases}\]
\[
\blw\left(t+1\right)=\calr\bls\left(t\right),\quad\forall t\geq0.
\]
In words, by eliminating an additional $1-\left(1-\mu_1 d\right)^2$ fraction of the population, we are able to set $\bls(0)$ to be a scalar multiple of $\blts(0)$. This is possible due to Eq.~\eqref{eq:lip2}. From $t=1$ and onwards, $\bls(t)$ is set to stay along the same direction as $\blts(t)$. One can verify that the construction guarantees that 
\begin{equation}
\label{eq:lip3}
	\rreal{\blw,\bls} \geq \rreal{\bltw,\blts}, \quad \forall t \geq 1.
\end{equation}
Denoting by $R_{\blw}$ the discounted total reward of the sequence $\left\{\left(\bltw\left(t\right),\blts\left(t\right)\right)\right\}$, and recalling that $\rreal{\blw,\bls}=\ln\left(\frac{\nor{\calr \bls}}{\nor{\blw}}\right)$, we have
\begin{eqnarray*}
R_{w} & = & \sum_{t=0}^{\infty}\gamma^{t}\rreal{\blw\left(t\right),\bls\left(t\right)}\\
& = & \ln \left(\frac{\nor{\calr \blts\left(0\right)}}{\nor{\bltw(0)}}\cdot\frac{\nor{\bltw(0)}}{\nor{\blw(0)}}\left(1-\mu_2 d\right)^2\right) + \sum_{t=1}^{\infty}\gamma^{t}\rreal{\blw(t),\bls(t)} \\
 & \stackrel{(a)}{\geq}& \rreal{\bltw(0),\blts(0)} + \ln\left(\frac{\nor{\bltw(0)}}{\nor{\blw(0)}}\left(1-\mu_2 d\right)^2\right) + \sum_{t=1}^{\infty}\gamma^{t}\rreal{\bltw(t),\blts(t)} \\
 &\stackrel{(b)}{\geq}& R_{\bltw} + \ln\left(1-\mu d\right)^3, 
\end{eqnarray*}
where $\mu=\max\left\{\mu_1,\mu_2\right\}$. Steps $(a)$ and $(b)$ follow from Eq.~\eqref{eq:lip3} and Eq.~\eqref{eq:lip1}, respectively. Since $R_{\blw} \leq V^\gamma\left(\blw\right)$, and $R_{\bltw} = V^\gamma\left(\bltw\right)$ for some feasible sequence starting from $\bltw$, we have $V^\gamma\left(\blw\right) - V^\gamma\left(\bltw\right) \geq \ln\left(1-\mu d\right)^3$. Repeating the same arguments with $\blw$ and $\bltw$ exchanged, we will have proven Eq.\eqref{eq:lip0}, which in turn establishes the claim. 
\Halmos\endproof

Finally, we show that the Lipschitz continuity of $V^{\gamma}\left(\cdot\right)$ can be extended from $\calc$ to the entire simplex $\simp$, again via a coupling argument. Let $\blw,\bltw\in\simp$, such that $\nor{\blw-\bltw}$
is suitably small. Let $\left\{ \left(\bltw\left(t\right),\blts\left(t\right)\right)\right\} _{t\geq0}$
be any feasible sequence with $\bltw\left(0\right)=\bltw$, and denote by $R_{\bltw}$ its discounted total reward. Consider a sequence $\left\{ \left(\blw\left(t\right),\bls\left(t\right)\right)\right\} _{t\geq0}$,
generated by
\[
\blw\left(0\right)=\blw,
\]
\[
\bls(t) \in\begin{cases}
\left\{ p\left(\blw\left(0\right),\blts\left(0\right)\right)\right\} , & \, t=0,\\
\arg\max_{\bls\in\phis{\blw\left(t\right)}}V^{\gamma}\left(L\left(\bls\right)\right), & t\geq1,
\end{cases}
\]
\[
\blw\left(t+1\right)=\calr\bls\left(t\right),\quad\forall t\geq0.
\]
where $p\left(\cdot,\cdot\right)$ is the projection of $\bls$ onto
$\phis{\blw}$, $p\left(\blw,\bls\right)\bydef\arg\min_{\bltw\in\phis{\blw}}\nor{\bltw-\bls}_{2}.$
In other words, we set $\bls\left(0\right)$ to be a closest point
to $\blts\left(0\right)$ in $\phis{\blw\left(0\right)}$, and then
follow the optimal stationary policy prescribed by the value function
from $t=1$ and onward. Recall from Eq.~\eqref{eq:phisexp} that 
\[
\phis{\blw} = \left\{ \bls \in\simp: 0\leq\bls_{i}\leq\frac{1}{\min\left\{\norp{\blw}, \nor{\blw}\right\}}\blw_{i},\,\forall1\leq i\leq K\right\}. 
\]
This implies that there exists $k_1 >0$ such that whenever $\nor{\blw\left(0\right)-\bltw\left(0\right)}$ is sufficiently small,
\begin{equation}
\nor{\bltw\left(1\right)-\blw\left(1\right)}\leq k_{1}\nor{\blw\left(0\right)-\bltw\left(0\right)}.\label{eq:v1w1}
\end{equation}
Also, by definitions of the reward $R\left(\cdot, \cdot\right)$ and that of $\bls(0)$, there exists $k_2>0$, such that whenever $\nor{\blw\left(0\right)-\bltw\left(0\right)}$ is sufficiently small,
\begin{eqnarray}
\left|R\left(\bltw(0), \blts\left(0\right)\right)-R\left(\blw(0), \bls\left(0\right)\right)\right| \leq k_{2}\nor{\blw\left(0\right)-\bltw\left(0\right)}. \label{eq:v0w0R}
\end{eqnarray}
Denoting by $R_{\blw}$ the discounted total reward from the sequence $\left\{ \left(\blw\left(t\right),\bls\left(t\right)\right)\right\}$,
we have that
\begin{eqnarray*}
R_{\blw} & = & \sum_{t=0}^{\infty}\gamma^{t}R\left(\blw(t),\bls\left(t\right)\right)\\
 & = & R\left(\bltw(t),\blts\left(0\right)\right)+\sum_{t=1}^{\infty}\gamma^{t}R\left(\blw(t),\bls\left(t\right)\right)+\left(R\left(\blw(0), \bls\left(0\right)\right)-R\left(\bltw(0), \blts\left(0\right)\right)\right)\\
 & = & \left(R\left(\bltw(0),\blts\left(0\right)\right)+\gamma V^{\gamma}\left(\blw\left(1\right)\right)\right)+\left(R\left(\blw(0),\bls\left(0\right)\right)-R\left(\bltw(0), \blts\left(0\right)\right)\right)\\
 & \stackrel{(a)}{\geq} & \left(R\left(\bltw(0), \blts\left(0\right)\right)+\gamma V^{\gamma}\left(\bltw\left(1\right)\right)\right)-\gamma l_{1}\nor{\blw\left(1\right)-\bltw\left(1\right)}\\
 && + \left(R\left(\blw(0), \bls\left(0\right)\right)-R\left(\bltw(0), \blts\left(0\right)\right)\right)\\
 & \stackrel{\left(b\right)}{\geq} & \left(R\left(\bltw(0), \blts\left(0\right)\right)+\gamma V^{\gamma}\left(\bltw\left(1\right)\right)\right)-\gamma k_{1}l_{1}\nor{\blw\left(0\right)-\bltw\left(0\right)}-k_{2}\nor{\blw\left(0\right)-\bltw\left(0\right)}\\
 & \geq & V^{\gamma}\left(\bltw\left(0\right)\right)-\left(\gamma k_{1}l_{1}+k_{2}\right)\nor{\blw\left(0\right)-\bltw\left(0\right)},
\end{eqnarray*}
where $(a)$ follows on Lemma \ref{lem:lipsinC}, and
the fact that $\blw\left(1\right),\bltw\left(1\right)\in\calc$, and $(b)$ from Eqs.~\eqref{eq:v1w1} and \eqref{eq:v0w0R}.
Finally, we have
\begin{eqnarray*}
V^{\gamma}\left(\bltw\left(0\right)\right)-V^{\gamma}\left(\blw\left(0\right)\right) & \leq & V^{\gamma}\left(\bltw\left(0\right)\right)-R_{\blw}\\
 & \leq & V^{\gamma}\left(\bltw\left(0\right)\right)-\left(-\left(\gamma k_{1}l_{1}+k_{2}\right)\nor{\blw\left(0\right)-\bltw\left(0\right)}+V^{\gamma}\left(\bltw\left(0\right)\right)\right)\\
 & = & \left(\gamma k_{1}l_{1}+k_{2}\right)\nor{\blw\left(0\right)-\bltw\left(0\right)}.
\end{eqnarray*}
Setting $l=\gamma k_{1}l_{1}+k_{2},$ this completes the proof of
Proposition \ref{prop:val_lipcont}. 

\subsubsection{From Discounted to Average-reward Criteria}

We show in this subsection that the concavity (Proposition \ref{prop:val_convx-1}) and smoothness (Proposition \ref{prop:val_lipcont}) of the \emph{discounted} value functions $V^\gamma(\cdot)$ imply the concavity and continuity of the bias function for the \emph{average-reward} problem. Define the {\bf average reward} of a policy $\pi$ starting from state $\blw$ at time $0$ as $J_{\pi}\left(\blw\right)=\limsup_{N\rightarrow\infty}\frac{1}{N}\sum_{n=0}^{N-1}R\left(\blw(t),\bls(t)\right),$ where $\left\{ \left(\blw(t), \bls\left(t\right)\right)\right\} _{t\geq0}$ is the induced state sequence
by applying policy $\pi$, with $\blw(0)=\blw$. Denote by 
\[J^{*}\left(\blw\right)=\sup_{\pi\in\Pi}J_{\pi}\left(\blw\right), \]
the optimal reward achievable from the point $\blw$. Since all entries of the matrix $\calr$ are positive, it is not difficult to see that the optimal reward is identical over the simplex, i.e.\ $\alpha^{*}=J^{*}\left(\blw\right),\,\forall\blw\in\simp$. We have the following characterization of an optimal policy. The proof is similar to that of the classical finite state space average-reward problem \cite{Ross83}. We include the proof in Appendix \ref{app:lem:avgBell} for completeness. 
\begin{lemma}
\label{lem:avgBell}If there exists a bounded function $g:\simp\rightarrow\R$
and a constant $\alpha$ such that 
\begin{equation}
\alpha+g\left(\blw\right)=\max_{\bls\in\phis{\blw}}\left[R\left(\bls\right)+g\left(\proj{\calr\bls}\right)\right].\label{eq:avgBell}
\end{equation}
Then a stationary policy $\pi^{*}$ that achieves the above maximum for all $\blw\in \simp$ also yields the optimal average reward, i.e., $J_{\pi^{*}}\left(\blw\right)=\alpha^{*},\,\forall\blw\in\simp.$
\end{lemma}

The following proposition is the main result of this subsection. 

\begin{proposition}
\label{prop:adv_val_conv}Fix any $\bltw\in\simp$. There exists an
increasing sequence $\gamma_{n}\nearrow1$, and function $g:\simp\rightarrow\R$
such that the following hold:
\begin{enumerate}
\item The function $g\left(\cdot\right)$ is bounded, concave and continuous. 
\item It holds that  $\lim_{n\rightarrow\infty}\sup_{\blw\in\Delta_{K}}\left(\left(V^{\gamma_{n}}\left(\blw\right)-V^{\gamma_{n}}\left(\bltw\right)\right)-g\left(\blw\right)\right)=0.$
\item The function $g\left(\cdot\right)$ satisfies Eq.~\eqref{eq:avgBell} for some
constant $\alpha$.
\end{enumerate}
\end{proposition}
\proof{Proof.}
Fix any $\bltw\in\simp$. Define 
\[
g^{\gamma}\left(\blw\right)\bydef V^{\gamma}\left(\blw\right)-V^{\gamma}\left(\bltw\right).
\]
Since $g^{\gamma}\left(\bltw\right)=0$ for all $\gamma$, by Proposition
\ref{prop:val_lipcont}, $g^{\gamma}\left(\cdot\right)$ is bounded
and $l$-Lipschitz-continuous for all $\gamma\in\left(\frac{1}{2},1\right).$
By the Arzela-Ascolli theorem for compact metric spaces, there exists
a Lipschitz-continuous function $g\left(\cdot\right)$ defined on
$\simp,$ and a positive increasing sequence $\gamma_{n}\rightarrow1$
such that 
\begin{equation}
\lim_{n\rightarrow\infty}\sup_{\blw\in\simp}\left(g^{\gamma_{n}}\left(\blw\right)-g\left(\blw\right)\right)=0.\label{eq:g_uniconv}
\end{equation}
By Proposition \ref{prop:val_convx-1}, $g^{\gamma_{n}}$ is concave for all $n$, and hence the limiting function $g$ is also concave. It remains to show that $g$ satisfies Eq.~\eqref{eq:avgBell}.
To this end, note that by Eq.~\eqref{eq:disBell}, for all $\blw\in\simp$,
\begin{eqnarray}
  \left(1-\gamma_{n}\right)V^{\gamma_{n}}\left(\bltw\right)+g^{\gamma_{n}}\left(\blw\right) & = & V^{\gamma_{n}}\left(\blw\right)-\gamma_{n}V^{\gamma_{n}}\left(\bltw\right)\nonumber \\
 & = & \max_{\bls\in\phis{\blw}}\left[R\left(\blw,\bls\right)+\gamma_n V^{\gamma_n}\left(\proj{\calr \bls}\right)\right] - \gamma_nV^{\gamma_n}(\bltw) \nonumber\\
 & = & \max_{\bls\in\phis{\blw}}\left[R\left(\blw,\bls\right)+\gamma_{n}g^{\gamma_{n}}\left(\proj{\calr\bls}\right)\right].\label{eq:gconv1}
\end{eqnarray}
Because $R$ is bounded on $\simp,$ we have that, for some $B>0$, 
$\left|V^{\gamma_{n}}\left(\blw\right)\right|$ is bounded over $\simp$ by $\frac{B}{1-\gamma_{n}}$. Hence, there exists a subsequence $\left\{ \gamma_{n_{k}}\right\} \subset\left\{ \gamma_{n}\right\} $
such that 
\begin{equation}
\lim_{k\rightarrow\infty}\left(1-\gamma_{n_{k}}\right)V^{\gamma_{n_{k}}}\left(\bltw\right)=\alpha,\label{eq:gconv2}
\end{equation}
for some $\alpha\in\R.$ Combining Eqs.~\eqref{eq:gconv1} and \eqref{eq:gconv2},
we have that
\begin{eqnarray*}
\alpha+g\left(\blw\right)&=&\lim_{k\rightarrow\infty}\left(1-\gamma_{n_{k}}\right)V^{\gamma_{n_{k}}}\left(\bltw\right)+g^{\gamma_{n_{k}}}\left(\blw\right)\\
 & = & \lim_{k\rightarrow\infty}\max_{\bls\in\phis{\blw}}\left[R\left(\blw,\bls\right)+\gamma_{n_{k}}g^{\gamma_{n_{k}}}\left(\proj{\calr\bls}\right)\right]\\
 & = & \max_{\bls\in\phis{\blw}}\left[R\left(\blw, \bls\right)+g\left(\proj{\calr\bls}\right)\right],
\end{eqnarray*}
where the last equality is based on the fact that $g^{\gamma_{n_{k}}}$
converges to $g$ uniformly over $\simp$ as $k\rightarrow\infty$
(Eq.~\eqref{eq:g_uniconv}). This completes the proof.
\Halmos\endproof

\subsubsection{Proof of Theorem \ref{thm:detcont}}
We are now ready to the establish the existence of an optimal population mixture. 
\proof{Proof.} {\bf (Theorem \ref{thm:detcont})}
Consider the SIM dynamic system. Denote by $L\left(\cdot\right)$
the map from the action, $\bls\left(t\right)$, to the next state, $\blw\left(t+1\right)$, 
\[
L\left(\bls\right)\bydef\proj{\calr\bls},\quad\forall\bls\in\simp,
\]
and by $H\left(\blw\right)$ the set of states in $\simp$ that are
reachable from $\blw$ in one step: $H\left(\blw\right)=\bigcup_{\bls\in\phis{\blw}}L\left(\bls\right)$.
Let $\Gamma\left(\cdot\right)$ be the set-valued map defined by 
\[
\gamap{\blw}\bydef \arg\max_{\bltw\in H\left(\blw\right)}g\left(\bltw\right),\quad\forall\blw\in\simp.
\]
By Lemma \ref{lem:avgBell}, it suffices to show that the
map $\gamap{\cdot}$ admits a fixed point, in the sense that $\blw\in\gamap{\blw}$, for some $\blw\in\simp$. To this end, we will invoke Kakutani's fixed point theorem \cite{Kak41}.\footnote{See Definition \ref{def:contSet} in Appendix \ref{app:techprelim} for different notions of continuity for set-valued functions.}
\begin{lemma}
\textbf{\label{lem:Kakutani}(Kakutani's Fixed-Point Theorem) }Let
$S$ be a non-empty, compact and convex subset of $\R^{K}.$ Let $\varphi:S\rightarrow2^{S}$
be an upper semicontinuous set-valued function on $S$, with the property
that $\varphi\left(x\right)$ is non-empty, closed and convex for
all $x\in S$. Then $\varphi$ admits a fixed point, i.e. $x\in \phi(x)$ for some $x\in S$. 
\end{lemma}
Clearly, the set $\simp$ is a non-empty, compact and convex subset
of $\R^{K}$. Using Lemma \ref{lem:Kakutani}, our proof will be completed
by showing the following. 
\begin{lemma}
\label{lem:gammacont}The map $\gamap{\cdot}$ is upper-semicontinuous
on $\simp$, and $\gamap{\blw}$ is non-empty, closed and convex for
all $\blw\in\simp.$
\end{lemma}

\proof{Proof.} {\bf (Lemma \ref{lem:gammacont})}
We first show the upper semicontinuity of $\gamap{\cdot}$. We will
use the following Maximum theorem (cf. \cite{Ok07}).
\begin{lemma}
\textbf{\label{lem:maximum} }(\textbf{Berge}'s \textbf{Maximum Theorem)}
Let $X$ and $\Sigma$ be metric spaces, $f:X\times\Sigma\rightarrow\R$
be a function jointly continuous in its two arguments, and $C:\Sigma\rightarrow2^{X}$
be a compact-valued map. For $\sigma\in\Sigma$, let
\[
C^{*}\left(\sigma\right)\bydef\arg\max_{x\in C\left(\sigma\right)}f\left(x,\sigma\right).
\]
If $C\left(\cdot\right)$ is continuous at some $\sigma$, then $C^{*}\left(\cdot\right)$
is non-empty, compact-valued, and upper semicontinuous at $\sigma.$
\end{lemma}
It can be checked that both mappings $L\left(\cdot\right)$ and $\phis{\cdot}$
are continuous. Hence, the composed map $H\left(\cdot\right)=L\cdot\phis{\cdot}$
is also continuous. By Proposition \ref{prop:adv_val_conv}, the bias function
$g\left(\cdot\right)$ is continuous on $\simp$. The upper semicontinuity
of $\gamap{\cdot}$ then follows from Lemma \ref{lem:maximum} by
letting $X=\Sigma=\simp$, $f\left(x,\sigma\right)=g\left(x\right),$
$C\left(\sigma\right)=H\left(\sigma\right)$, and $C^{*}\left(\sigma\right)=\gamap{\sigma}.$ 

It remains to be shown that $\gamap{\cdot}$ is convex-valued at all
$\blw\in\simp$. We first claim that the set $H\left(\blw\right)$
is convex for all $\blw\in\simp.$ Since it can be easily checked
from the definition that $\phis{\blw}$ is convex for all $\blw\in\simp$,
it suffices to show that for all $\blv,\blu\in\simp$ and all $a \in [0,1]$,
there exists $b \in [0,1]$ such that 
\[
aL\left(\blu\right)+\left(1-a\right)L\left(\blv\right)=L\left(b\blu+(1-b)\blv\right),
\]
which would then imply that the image of a convex subset of $\simp$
under $L\left(\cdot\right)$ remains convex. By the definition of
$L\left(\cdot\right)$, we have that 
\begin{eqnarray*}
  L\left(b\blu+(1-b)\blv\right) &=& \proj{\calr\left(b\blu+(1-b)\blv\right)} = \proj{b\calr\blu+(1-b)\calr\blv}\\
 & = & \frac{b\calr\blu+(1-b)\calr\blv}{\nor{b\calr\blu+(1-b)\calr\blv}}\\
 & \stackrel{(a)}{=} & \frac{b\nor{\calr\blu}L\left(\blu\right)+\left(1-b\right)\nor{\calr\blv}L\left(\blv\right)}{\nor{b\calr\blu+(1-b)\calr\blv}}\\
 & = & f(b)L\left(\blu\right)+(1-f(b))L\left(\blv\right),
\end{eqnarray*}
where $f(b)\bydef\frac{b\nor{\calr\blu}}{\nor{b\calr\blu+(1-b)\calr\blv}}$, and $(a)$ follows from the fact that $L(\blu)=\blu$ for all $\blu\in\simp$. Since all entries of $\calr$ are positive, $f$ is continuous over the interval $[0,1]$, $f(0)=0$, and $f(1)=1$. By the intermediate value theorem, there exists some $b$ such that $f(b)=a$. This completes the proof of Lemma \ref{lem:gammacont}.
\Halmos\endproof

The proof of Theorem \ref{thm:detcont} is completed by combining
Lemmas \ref{lem:Kakutani} and \ref{lem:gammacont}. \Halmos\endproof

\section{Proof of Theorem \ref{thm:k2}}
\label{sec:k2pf}

In this section, we prove Theorem \ref{thm:k2}, which gives an explicit formula for the fixed point corresponding to the optimal stationary policy when $K=2$. The main argument is that of proof by contradiction: for any other fixed point, we will construct a superior policy that provides a strictly greater average reward, thus invalidating its optimality. We note that the proof for Theorem \ref{thm:k2} heavily relies the fact that $K=2$: the per-stage reward function, $R\left(\blw,\cdot\right)$, is \emph{monotone} over the simplex $\simp$, which in this case is just a one-dimensional line segment. Unfortunately, when $K\geq 3$, the geometry of $\simp$ becomes non-trivial, and the evolutions of the state sequence $\blw(t)$ are more complex. As a result, our current proof-by-contradiction arguments do not appear to generalize to higher dimensions. 

\proof{Proof.} {\bf (Theorem \ref{thm:k2})}
We first observe that the result trivially holds when $\nor{\calr_{(\cdot,1)}}=\nor{\calr_{(\cdot,2)}}$ (in which case all points on $\simp$ are optimal), or when $\calr_{(\cdot,1)}=a \calr_{(\cdot,2)}$ for some $a>0$. Therefore, without loss of generality, for the rest of the proof we will assume that $\nor{\calr_{(\cdot,1)}}>\nor{\calr_{(\cdot,2)}}$ and that $\calr_{(\cdot,1)}$ and $\calr_{(\cdot,2)}$ are linearly independent. 

Consider the dynamic system SIM. Since $K=2$, the action set $\phis{\blw}$ at state $\blw\in\simp$ is a one-dimensional line segment, and we can re-write it as\footnote{Notation: $\blw_i$ represents the $i$th coordinate of the vector $\blw$.}
\[
\phis{\blw}=\left\{ \bls\in\simp:l_{\beta}\left(\blw_{1}\right)\leq\bls_{1}\leq h_{\beta}\left(\blw_{1}\right)\right\} ,
\]
where $h_{\beta}(x)$ and $l_{\beta}(x)$ represent the lower and higher end of the action set, respectively, when $\blw_1=x$. It is easy to verify that
\begin{eqnarray}
h_{\beta}(x) & = & \min\left\{ \frac{x}{\beta},1\right\} ,\label{eq:h}\\
l_{\beta}(x) & = & \max\left\{ \frac{x}{\beta}-\left(\frac{1}{\beta}-1\right),0\right\} .\label{eq:l}
\end{eqnarray}
We have the following characterization of $h_{\beta}(\cdot)$ and $l_{\beta}(\cdot)$, which follows immediately by definition. 
\begin{lemma}
\label{lem:hlmon}$h_{\beta}\left(\cdot\right)$ and $l_{\beta}\left(\cdot\right)$
are continuous and non-decreasing on the interval $[0,1]$.\end{lemma}
%
%


The next lemma shows that the sub-population chosen in a stationary optimal policy must lie on the boundary of the action set. The proof is given in Appendix \ref{app:lem:boundary}. 
\begin{lemma}
\label{lem:boundary}For all $K\geq2$, let $\pi^{*}$ be an optimal
stationary policy, and let $\blw$ be an optimal fixed-point under $\pi^{*}$,
in the sense of Eq.~\eqref{eq:fixptdef}. We have that $\pi\left(\blw \right)\in\partial\left(\phis{\blw}\right)$, where $\partial\left(X\right)$denotes the boundary of the set $X$. \end{lemma}
Let $M(x)$ be the first coordinate of $\blw(t+1)$ when the sub-population $\bls(t) = (x,1-x)^\top$, i.e.,
\[M\left(x\right)\bydef\left(\proj{\calr\left(x,1-x\right)^{\top}}\right)_{1}.\]
We have the following two lemmas. The proofs are elementary and is omitted. 

\begin{lemma}
\label{lem:M-mon} The function $M(x)$ is continuously differentiable over $(0,1)$, and 
\[
 \left(\proj{\rcol 1}-\proj{\rcol 2}\right)_{1} \cdot \left(\frac{d}{dx}M(x)\right)>0,\quad\forall 0\leq x\leq1.
\]
\end{lemma}

\begin{lemma}
\label{lem:costmon} Consider $R\left(\blw, \bls\right)$, the per-stage reward function for SIM. For all $\blw,\bls \in (0,1)$, we have that $\frac{\partial}{\partial\bls_{1}}R\left(\blw,\bls\right)= \frac{\partial}{\partial\bls_{1}}\ln \left(\nor{\calr\bls}\right) <0$. 
\end{lemma}

We are now ready to characterize the form of the optimal stationary policy, $\pi^*$. We consider the following two cases, depending on the value of $\left(\proj{\rcol 1}-\proj{\rcol 2}\right)_{1}$.\footnote{Note that due to the assumption on the linear independence between
$\rcol 1$ and $\rcol 2$, we have that $\left(\proj{\rcol 1}-\proj{\rcol 2}\right)_{1}\neq0$.}

%

\begin{figure}[ht]
\begin{minipage}[b]{0.45\linewidth}
\centering \includegraphics[scale=0.45]{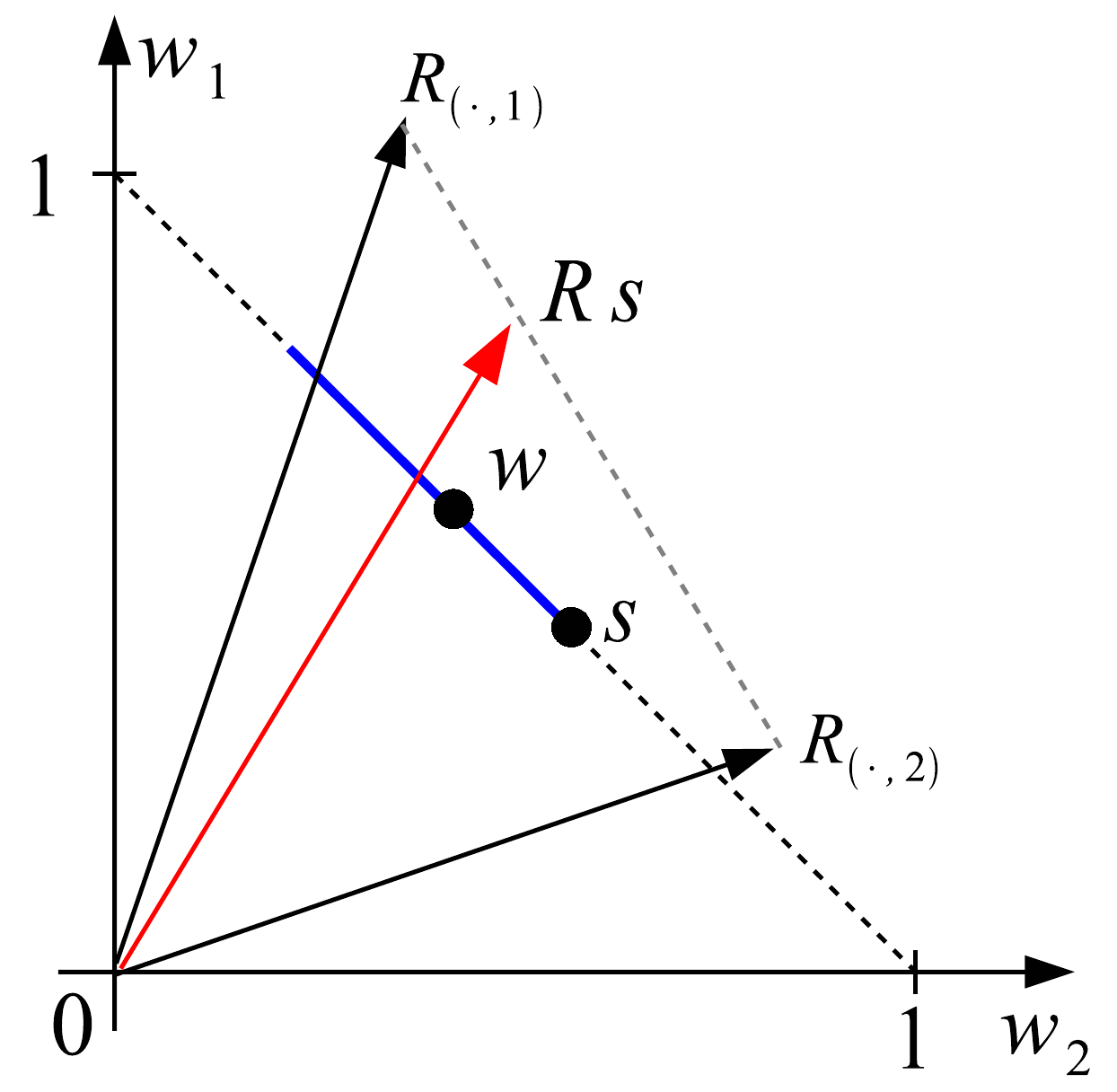} \caption{ $\left(\proj{\rcol 1}-\proj{\rcol 2}\right)_{1}>0$.}
\label{fig:2dcase2} 
\end{minipage}
\hspace{0.5cm}
\begin{minipage}[b]{0.45\linewidth}
\centering 
\includegraphics[scale=0.45]{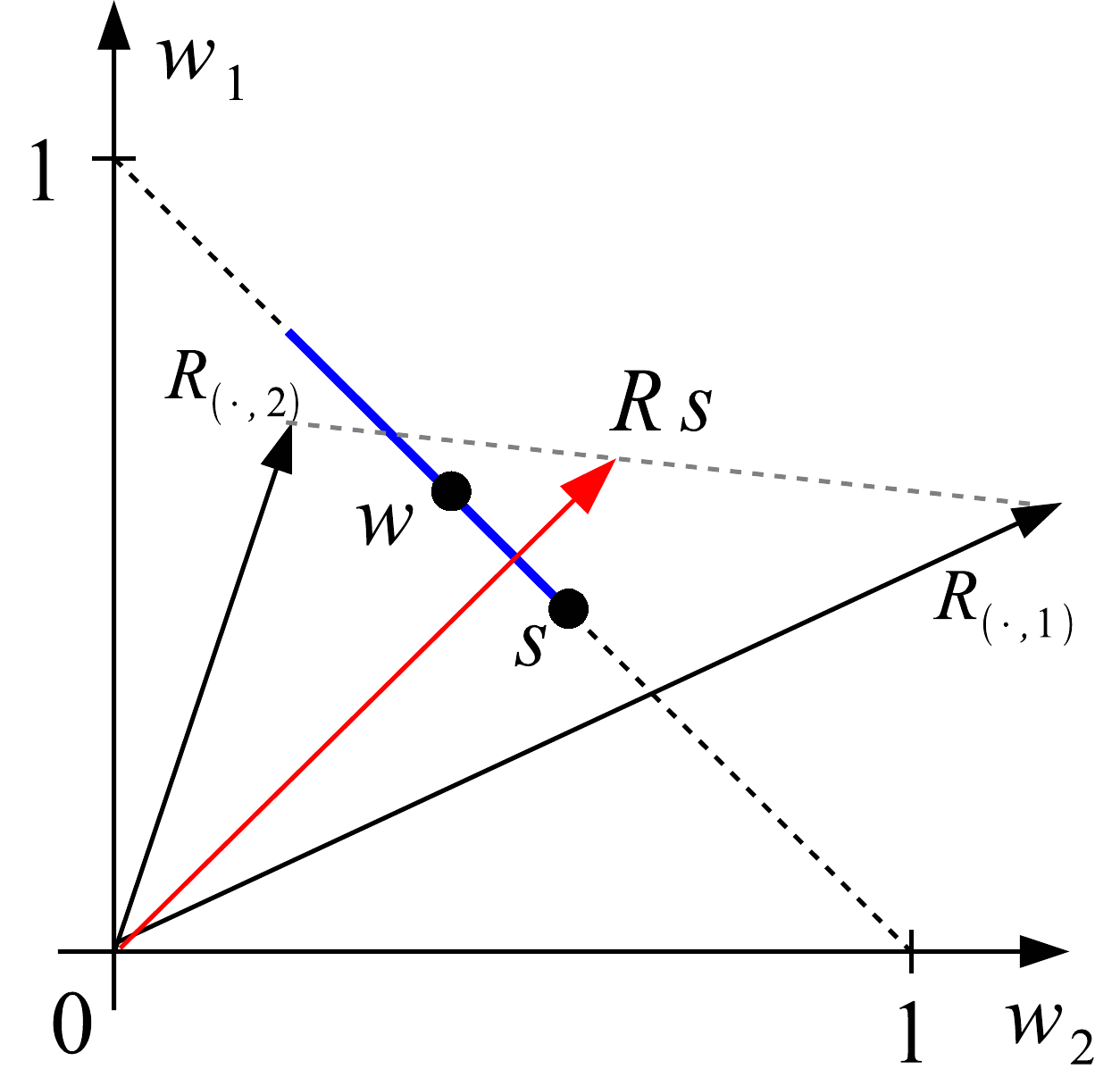} 
\caption{$\left(\proj{\rcol 1}-\proj{\rcol 2}\right)_{1}<0$.}
\label{fig:2dcase3} 
\end{minipage}
\end{figure}
\vspace{-10pt}

\begin{claim}
\label{clm:h}If $\left(\proj{\rcol 1}-\proj{\rcol 2}\right)_{1}>0$ (Figure \ref{fig:2dcase2}),
then $\left(\pi^{*}\left(\blw^{*}\right)\right)_{1}=h_{\beta}\left(\blw_{1}^{*}\right).$\end{claim}
\proof{Proof.}
By Lemma \ref{lem:boundary}, we know that $\left(\pi^{*}\left(\blw^{*}\right)\right)_{1}\in\left\{ \lowpt{\blw_{1}^{*}},\highpt{\blw_{1}^{*}}\right\} .$
Suppose for contracdiction that $\left(\pi^{*}\left(\blw^{*}\right)\right)_{1}=l_{\beta}\left(\blw_{1}^{*}\right)$
for some $\blw^{*}\in\simp$. We show that there exists a fixed point
for some other policy $\hat{\pi}$ with a strictly greater reward.
It is easy to verify that if $\blw=(0,1)^\top$ or $(1,0)^\top$, then $\phis{\blw}=\left\{\blw\right\}$. Since $\calr_{\left(i,j\right)}>0$ for all $i,j,$, this implies that if $\blw^{*}$ is a fixed point, then $\blw_{1}^{*}\neq0\mbox{ or }1.$ Therefore, we will assume $0<\blw_{1}^{*}<1$. We have that 
\begin{eqnarray}
\blw^*_{1}-\lowpt{\blsw_{1}} & = & \blsw_{1}-\max\left\{ \frac{\blsw_{1}}{\beta}-\left(\frac{1}{\beta}-1\right),0\right\} = \min\left\{ \blsw_{1},\left(1-\blsw_{1}\right)\left(\frac{1}{\beta}-1\right)\right\} \nonumber \\
 & > & 0.\label{eq:wgap}
\end{eqnarray}
Let $\blts\in\simp$ be given by $\blts_{1}=\frac{\blw^*_{1}+\lowpt{\blsw_{1}}}{2}$, and let the map $H_{l}$ be defined by 
\begin{equation}
H_{l}\left(x\right)\bydef\lowpt{M\left(x\right)},\quad0\leq x\leq1.\label{eq:H_l}
\end{equation}
By Lemmas \ref{lem:hlmon} and \ref{lem:M-mon}, $H_{l}\left(x\right)$
is continuous and monotonically non-decreasing on the interval $\left[0,1\right].$
Therefore, we have that
\[
\blts_{1}\leq H_{l}\left(x\right)\leq1,\quad\forall \blts_{1}\leq x\leq1.
\]
By Brouwer's fixed point theorem, this implies that there exists $\blts^{*}\in\simp$
such that $\blts_{1}\leq\blts_{1}^{*}\leq1$ and $H_{l}\left(\blts_{1}^{*}\right)=\blts_{1}^{*}.$
By Eq.~\eqref{eq:H_l}, the point $\bltw^{*}=\proj{\calr\blts^{*}}$
is a fixed point for the policy $\hat{\pi}$ that always choose the point $\lowpt{\blw_1}$:
\[
\left(\hat{\pi}\left(\blw\right)\right)_{1}=\lowpt{\blw_{1}},\quad\blw\in\simp.
\]
By Eq.~\eqref{eq:wgap}, we have that 
\[
\blts_{1}^{*}=\left(\hat{\pi}\left(\blv^{*}\right)\right)_{1}>\left(\pi\left(\blw^{*}\right)\right)_{1}=\lowpt{\blw_{1}^{*}}.
\]
By Lemma \ref{lem:costmon}, we have that the fixed point $\bltw^{*}$ under
$\hat{\pi}$ provides a strictly greater average reward compared to
that of $\blw^{*}$under $\pi^{*}$, which contradicts with the optimality
of $\blw^{*}$. This completes the proof of the claim. \Halmos\endproof

\begin{claim}
\label{clm:l}If $\left(\proj{\rcol 1}-\proj{\rcol 2}\right)_{1}<0$ (Figure \ref{fig:2dcase3}),
then $\left(\pi^{*}\left(\blw^{*}\right)\right)_{1}=\highpt{\blw_{1}^{*}}.$
\end{claim}

\proof{Proof.} Similar to
the proof of Claim \ref{clm:h}, assume for contradiction that $\left(\pi^{*}\left(\blw^{*}\right)\right)_{1}=\lowpt{\blw_{1}^{*}}.$
With a similar argument to that of Eq.~\eqref{eq:wgap}, we have that
\begin{equation}
\highpt{\blw_{1}}>\lowpt{\blw_{1}},\quad\forall0<\blw_{1}<1.\label{eq:highlow}
\end{equation}
Since $\left(\proj{\rcol 1}-\proj{\rcol 2}\right)_{1}<0$,
by Lemma \ref{lem:M-mon}, we have that $\frac{d}{dx}M(x)<0.$ By
Eq.~\eqref{eq:highlow}, there exists $\blts\in\simp$ such that
$\highpt{\blsw_{1}}>\blts_{1}>\lowpt{\blw_{1}^{*}}$ and 
\[
\highpt{\proj{\calr\blts}}>\lowpt{\blsw_{1}}>\lowpt{\proj{\calr\blts}},
\]
or in other words, $\lowpt{\blsw_{1}}\in\phis{\proj{\calr\blts}}$. Consider any stationary policy $\hat{\pi}$ that satisfies
\[
\left(\hat{\pi}\left(\blw\right)\right)_{1}=\blts_{1}\mbox{ and }\left(\hat{\pi}\left(\proj{\calr\blts}\right)\right)_{1}=\lowpt{\blw_{1}^{*}}.
\]
 Note that because of the assumption that $\blw^{*}$is a fixed point
for $\pi^{*},$ we have that 
\[
\proj{\calr\left(\lowpt{\blsw},1-\lowpt{\blsw}\right)^{\top}}=\blsw.
\]
 Starting from the point $\blw\left(0\right)=\blw^{*}$, it is not
difficult to see that the policy $\hat{\pi}$ produces the following
sequence of $\tilde{\blw}\left(t\right)$:
\[
\tilde{\blw}\left(t\right)=\begin{cases}
\left(\lowpt{\blsw},1-\lowpt{\blsw}\right)^{\top}, & t\mbox{ is even,}\\
\blts, & t\mbox{ is odd.}
\end{cases}
\]
By Lemma \ref{lem:costmon}, $\hat{\pi}$ produces a sequence of states
that yield a strictly higher average reward than the fixed point $\blw^{*}$.
This proves the claim. 
\Halmos\endproof
Theorem \ref{thm:k2} follows directly by combining Claims \ref{clm:h}
and \ref{clm:l}. \Halmos\endproof

%
%

\section{Proof of Theorem \ref{thm:probablistic}}
\label{sec:probpf}

We prove Theorem \ref{thm:probablistic} in this section, which shows that the results derived under the deterministic setting 
carry also into the stochastic case, where the set of offspring for each reproductive individual is random. Analogous to the case of classical branching processes, the optimal growth rate $\alpha^*$ associated with the reproduction matrix turns out to be an important qualitative property of the corresponding stochastic branching process, in the sense that there is a positive probability for the population to grow unboundedly if $\alpha^*>0$, while the process becomes extinct with probability one if $\alpha^*<0$. 

We provide some intuition as for why the growth rate properties of the deterministic model should remain valid for the stochastic case.  
The deterministic branching process model makes two critical assumptions,
namely, that the population is infinitely divisible, and that the
number of offspring produced by each individual is deterministic. Therefore, when $\alpha^*>0$, the deterministic approximation becomes more accurate in the regime of \emph{large populations}, so that the set of possible sub-population
(scaled to the unit simplex) can be approximated by its continuous
counterpart, and the stochastic fluctuations from reproduction becomes sufficiently small compared to the size of the population. As a result, the numbers of offspring produced in each type in each round concentrate around their
expected values, due to the Law of Large Numbers. On the other hand, if $\alpha^*<0$, we will then study the evolution of the \emph{expected} population process, $\left\{\E{\buz(t)}\right\}_{t\geq0}$, and show that the sequence $\left\{\E{\buz(t)}\right\}_{t\geq 0}$ always constitutes a feasible deterministic branching process in REAL. Since $\alpha^*<0$, we know that $\E{\buz(t)}$ must approach zero exponentially fast. 

%
%

\proof{Proof.} {\bf (Theorem \ref{thm:probablistic})} We will use the growth factor $\kappa^* =  e^{\alpha^*}$ instead of $\alpha^*$ in the proof, since it is more natural for mathematical manipulations in this context. We consider two cases, $\kappa^*<1$ and  $\kappa^*>1$, corresponding to $\alpha^*<0$ and  $\alpha^*>0$, respectively. 

Case $1$:  $\kappa^*<1$. The idea is to use the fact that the
sequence of \emph{expected} populations and sub-populations induced by any policy forms
a feasible deterministic branching process in REAL. We can then use the fact
that $\kappa^*<1$ to argue that $\lim_{t\rightarrow\infty}\E{\nor{\buz\left(t\right)}}\leq c\cdot\lim_{t\rightarrow\infty}\left(\kappa^*\right)^{t}=0$.
To make this argument concrete, fixing any policy $\pi_{t}\left(\cdot\right)$, 
we have that $ \E{\buz\left(t+1\right)}=\calr\E{\pi_{t}\left(\buz\left(t\right)\right)}$. 
Since $\nor{\pi_{t}\left(\buz\left(t\right)\right)}\leq\beta\nor{\buz\left(t\right)}$,
and $0\leq\left(\pi_{t}\left(\buz\left(t\right)\right)\right)_{i}\leq\buz_{i}\left(t\right)$, $\forall 1\leq i\leq K$,
we have that 
\begin{equation}
\E{\pi_{t}\left(\buz\left(t\right)\right)}\in\phir{\E{\buz\left(t\right)}}.\label{eq:expfeas}
\end{equation}
Let $\blw\left(t\right)=\E{\buz\left(t\right)}$ and $\bls\left(t\right)=\E{\pi_{t}\left(\buz\left(t\right)\right)}$.
Eq.~\eqref{eq:expfeas} implies that $\left\{ \left(\blw\left(t\right),\bls\left(t\right)\right)\right\} _{t\geq0}$
is a feasible sequence for the corresponding deterministic
branching process in REAL. Since $\kappa^*$is the optimal growth factor for
the deterministic process, we have that 
\begin{equation}
\lim_{t\rightarrow0}\nor{\E{\buz\left(t\right)}}\leq c\cdot\lim_{t\rightarrow\infty}\left(\kappa^*\right)^{t}=0.\label{eq:meanto0}
\end{equation}
By Eq.~\eqref{eq:meanto0}, for any finite initial condition,
there exists some decreasing function $g(t)$, $\lim_{t\rightarrow \ity}g(t)=0$, such that
\begin{equation}
	\pb\left(\buz\left(t\right)=0 \mbox{ for some } t \geq 0 \right) \geq \liminf_{t\rightarrow \ity} \pb\left( \buz(t) =0\right) \stackrel{(a)}{\geq} 1- \limsup_{t\rightarrow0} \nor{\E{\buz\left(t\right)}}  = 1,
\end{equation}
where $(a)$ follows from Markov's inequality by noting that $\nor{\buz(t)}$ must be integer-valued:
\[\pb\left( \nor{\buz(t)} =0\right)  = 1- \pb\left( \nor{\buz(t)} \geq  1 \right)\leq 1- 1\E{\nor{\buz(t)}}  = 1- \nor{\E{\buz(t)}}.\]

Case $2$: $\kappa^*>1$. We will construct a non-stationary policy
$\pi$ that achieves explosion with positive probability. Note that arguments involving the analysis of $\E{\buz(t)}$, similar to those above for the case $\alpha^*<1$, can be used to show that, with probability one, the population cannot explode at any rate strictly greater than $\kappa^*$. Therefore, we will focus on showing that the rate $\kappa^*$ is indeed achievable.

Let $\blx^{*}$ and $\bls^{*}$ be the optimal fixed-point and the corresponding sub-population, as was specified in Theorem \ref{thm:detcont}.
Starting with a sufficiently large population, the idea is to stay
along the direction of the optimal population mixture $\blx^{*}$,
by repeatedly finding a sub-population profile that is closest along
the direction of $\bls^{*}$. By the definition of $\bls^{*}$,
the resulting growth rate will be $\kappa^*$, provided that $\bus(t)$ stays closely along the direction of $\bls^*$ for all large $t$.

Let $\min_{1\leq i \leq K} \buz_{i}\left(0\right)$ be sufficiently large. In each time slot $t\geq0$, let%
\footnote{Here $\ceils{\blx}\bydef\left(\ceils{\blx_{1}},\ceils{\blx_{2}},\ldots,\ceils{\blx_{K}}\right)$.%
} 
\[
\bar{k}=\sup\left\{ k>0:\,\ceils{k\bls^{*}}\in\Phi\left(\buz\left(t\right)\right)\right\} ,
\]
where $\bls^{*}$ is the corresponding sub-population corresponding
to the optimal mixture $\blx^{*}$, and $\Phi\left(\cdot\right)$ was defined in Eq.~\eqref{eq:phiz}. Our \emph{policy} will be to choose the subpopulation
\[
\bus\left(t\right)=\ceils{\bar{k}\bls^{*}}, \quad  \forall t\geq0. 
\]
Let $A_{t}\left(\delta\right)$ be the event that 
\[
A_{t}\left(\delta\right)=\left\{ \buz_{i}\left(t+1\right)\geq\left(1-\delta\right)\left(\calr\bus\left(t\right)\right)_{i},\:\forall1\leq i\leq K\right\} ,
\]
and let $A=\cap_{t\geq0}A_{t}\left(\delta_{t}\right)$, with $\delta_{t}=\sqrt{\frac{ct}{k^{t}}}$
for positive some constants $k>1$ and $c>0$. We first show that
for all sample paths in $A$ (if $A$ is non-empty), the asymptotic
growth factor is no less than $\kappa^*$. To do so, we will lower
bound the sample paths in $A$, $\left\{ \left(\buz\left(t\right),\bus\left(t\right)\right)\right\} _{t\geq0}$,
by a sequence for REAL, $\left\{ \left(\blw\left(t\right),\bls\left(t\right)\right)\right\} _{t\geq0}$,
which achieves the growth factor $\kappa^*$. Let $\blw\left(0\right)=\buz\left(0\right)$, and
\[
\bls\left(t\right)=\left[\nor{\blw\left(t\right)}\left(\tilde{\beta}\left(t\right)-\frac{2}{\nor{\blw\left(t\right)}}\right)\left(1-\delta_{t}\right)\right]\bls^{*},
\]
where 
\[
\tilde{\beta}\left(t\right)=\begin{cases}
\max\left\{\beta \in [0,1]: \bls^* \in \phif{\blw(0),\beta}\right\} , & t=0,\\
1, & t\geq1,
\end{cases}
\]
and $\phif{\blw,\beta}$ was defined in Eq.~\eqref{eq:phif2}. The $-\frac{2}{\nor{\blw\left(t\right)}}$ adjustment is made
to take into account the quantization effect in the stochastic sample paths (i.e., $\buz_i(t)$ being integer-valued). One can show, by the construction of $\left\{ \blw\left(t\right)\right\} _{t\geq0}$,
that given a sufficiently large
initial population $\buz\left(0\right)$, 
\begin{equation}
	\bus_{i}\left(t\right)  \geq  \bls_{i}\left(t\right), \mbox{ and } \buz_{i}\left(t\right) \geq  \blw_{i}\left(t\right),
\end{equation}
for all sample paths in $A$, $i\in \left\{1,\ldots,K\right\}$, and $t\geq0$. The growth rate associated with $\left\{ \blw\left(t\right)\right\} _{t\geq0}$
is given by 
\[
\limsup_{t\rightarrow\infty}\frac{1}{t}\ln\nor{\blw\left(t\right)}=\kappa^*+\lim_{t\rightarrow \ity}\frac{1}{t}\sum_{s=0}^{t-1} \left[\ln\left(1-\delta_{s}\right)+ \ln\frac{1}{\beta}\left(\tilde{\beta}\left(s\right)-\frac{2}{\nor{\blw\left(s\right)}}\right)\right]=\kappa^*,
\]
where the last equality comes from the fact that $\lim_{t \to \infty}\delta_t= 0$ and $\lim_{t\to \infty} \nor{\blw(t)} = \ity$. This shows that for all sample paths in $A$, $\alpha=\limsup_{t\rightarrow\infty}\frac{1}{t}\ln\nor{\buz\left(t\right)}\geq\kappa^*$ .

Finally, we show that $\pb\left(A\right)>0$, for appropriately
chosen values of $\delta_t$. Given a sufficiently large initial population $\buz(0)$, there exists $\delta'>0$ and $k'>1$, such that if $\delta_t\leq \delta'$ for all large $t$, then for all sample paths in $A$, 
\begin{equation}
\label{eq:minRS}
\min_{1\leq i\leq K}\left(\calr\bus\left(t\right)\right)_{i}\geq {k'}^{t}, \forall t \geq 0. 
\end{equation} 
By an argument involving the Hoeffding's inequality (cf. Lemma \ref{lem:hoeffding} in Appendix \ref{app:techprelim}), and a union bound across all types, one can show that there exists
some constant $a,b>0$ so that
\begin{equation}
\pb\left(A_{t}\left(\delta\right)\right)\geq1-a\exp\left(-b\cdot\delta^{2}\min_{i}\left(\calr\bus\left(t\right)\right)_{i}\right),\label{eq:chernoff}
\end{equation}
for all $t\geq0$. We have that 
\vspace{-10pt}
\begin{align*}
\pb\left(A\right)  = & \pb\left(\cap_{t\geq0}A_{t}\left(\delta_{t}\right)\right) = 1 - \pb\left(\cup_{t\geq0}\overline{A_{t}\left(\delta_{t}\right)}\right) \\
  \stackrel{\left(a\right)}{\geq} & 1- \left(1-\pb\left(A_0\left(\delta_0\right)\right)\right)+\sum_{t=1}^{\infty}\left(1-\pb\left(A_{t}\left(\delta_{t}\right)\Big|\cap_{0\leq s\leq t-1}A_{s}\left(\delta_{s}\right)\right)\right)\\
  \stackrel{\left(b\right)}{\geq} & 1-\sum_{t=0}^{\infty}a\exp\left(-b\cdot\delta_{t}^{2}{k'}^{t}\right)=  1-\frac{a}{1-e^{-bc}},
\end{align*}
where ineuqality $\left(a\right)$ follows from a union bound argument, and
inequality $\left(b\right)$ from Eqs.~\eqref{eq:minRS} and \eqref{eq:chernoff}.
The last equality from the fact that $\delta_{t}=\sqrt{\frac{ct}{k^{t}}}$, with $k=k'$. We finish by choosing $c$ to be sufficiently
small so that $\pb\left(A\right)\geq1-\frac{a}{1-e^{-bc}}>0$. This completes the proof of Theorem \ref{thm:probablistic}. \Halmos\endproof

\section{Conclusions and Future Work}
\label{sec:concl}
We have shown in this paper that, in a class of multi-type branching processes with linear resource constraints, the optimal growth rate of the population can be achieved by a single appropriately chosen population mixture. The result holds both for deterministic and stochastic branching processes. For the special case of a two-type population where both types generate the same revenue, the optimal population mixture is obtained in closed form. 

There are many open questions worth exploring for the future. On the computational side, while the optimal mixture $\blx^*$ can be found by naively searching through all sub-populations $\bls$ that satisfy $\bls \in \proj{\phir{\left(\calr\bls\right)}}$, the complexity of doing so grows exponentially with $K$. We have not been able to establish an analog of the characterization in Theorem \ref{thm:k2} for $K\geq 3$ and non-symmetric revenue-per-individual, but results along these lines are clearly important for making the computation of the optimal mixture feasible in higher dimensions. Another interesting direction is to consider the ``dual'' problem of minimizing the growth rate, subject to a similar linear constraint (e.g., one can \emph{at most} eliminate $\beta$ fraction of the current population), with potential applications in vaccination \cite{FKG03}, queuing theory \cite{Res93}, etc. 

On the practical end, an important issue is that the values of the reproduction matrix $\calr$ can be difficult to obtain. For instance, as pointed out in \cite{DH06}, the rates at which quiescent cells become active are often difficult to measure precisely. This, along with other obstacles associated with the real-world manifestation of the current framework, such as $\calr$ being partially observed or time-varying, should bring a multitude of new challenges.

\newpage
\normalsize
\begin{APPENDIX}{}

\section{Technical Preliminaries}
\label{app:techprelim}

\begin{definition}
\label{def:minkowski}
The Minkowski addition of two sets $A,B\subset\R^{K}$ is defined as 
\[
A+B\bydef\left\{ \blx\in\R^{K}:\blx=\mathbf{a}+\mathbf{b},\mbox{ for some }\mathbf{a}\in A,\mbox{ }\mathbf{b}\in B\right\},
\]
The multiplication of $A\subset \R^K$ by a scalar $b \in \R$ is defined as
\[bA \bydef\left\{ \blx\in\R^{K}:\blx=b \mathbf{a},\mbox{ for some }\mathbf{a}\in A \right \}.\]
\end{definition}

\begin{lemma}
\label{lem:localglobal} A function $f\left(\cdot\right)$ that is $\epsilon$-locally $l$-Lipschitz continuous is also $l$-Lipschitz continuous.\end{lemma}

\proof{Proof.} We prove the lemma in the case of an $L_1$ norm, but the claim also holds for other norms as well. Fix any $\blx,\bly\in\R^{K}$ and suppose that $f$$\left(\cdot\right)$
is $\epsilon$-locally $l$-Lipschitz continuous. Let $N=\left\lceil \frac{\nor{\blx-\bly}}{\epsilon}\right\rceil $,
and let 
\[
\blu\left(i\right)=\blx+i\frac{\bly-\blx}{N},\quad1\leq i\leq N-1.
\]
We have that 
\begin{eqnarray*}
 &  & \left|f\left(\blx\right)-f\left(\bly\right)\right|\\
 & \stackrel{}{\leq} & \left|f\left(\blu\left(1\right)\right)-f\left(\blx\right)\right|+\left(\sum_{i=2}^{N-1}\left|f\left(\blu\left(i\right)\right)-f\left(\blu\left(i-1\right)\right)\right|\right)\\
 &  & +\left|f\left(\bly\right)-f\left(\blu\left(N-1\right)\right)\right|\\
 & \stackrel{(a)}{\leq} & l\left(\nor{\blu\left(1\right)-\blx}+\left(\sum_{i=2}^{N-1}\nor{\blu\left(i\right)-\blu\left(i-1\right)}\right)+\nor{\bly-\blu\left(N-1\right)}\right)\\
 & \stackrel{(b)}{=} & l\nor{\blx-\bly},
\end{eqnarray*}
where $(a)$ follows from the $\epsilon$-locally $l$-Lipschitz
continuity of $f\left(\cdot\right)$, and $\left(b\right)$
from the use of an $l_{1}$ norm, and the fact that the points $\left\{ \blu\left(1\right),\blu\left(2\right),\ldots,\blu\left(N-1\right)\right\} $
lie on the line segment joining $\blx$ and $\bly$. 
\Halmos\endproof

\begin{definition} \label{def:contSet}
A set-valued function $\Gamma:X\rightarrow2^{X}$ is said to be upper
semicontinuous, lower semicontinuous, and continuous at $x,$ if for
any sequence $\left\{ x_{n}\right\} _{n\geq0}\subset X$ such that
$\lim_{n\rightarrow\infty}x_{n}=x$, we have that 
\begin{eqnarray*}
\Gamma\left(x\right) & \supset & \limsup_{n\rightarrow\infty}\Gamma\left(x_{n}\right),\\
\Gamma\left(x\right) & \subset & \liminf_{n\rightarrow\infty}\Gamma\left(x_{n}\right),\mbox{ }\\
\Gamma\left(x\right) & = & \left(\limsup_{n\rightarrow\infty}\Gamma\left(x_{n}\right)\right)\cap\left(\liminf_{n\rightarrow\infty}\Gamma\left(x_{n}\right)\right),
\end{eqnarray*}
respectively, where $\limsup_{n\rightarrow\infty}X_{n} \bydef \bigcap_{n=1}^\infty \bigcup _{m=n}^\infty X_m$, and $
\liminf_{n\rightarrow\infty}X_{n} \bydef \bigcup_{n=1}^\infty \bigcap _{m=n}^\infty X_m$.
%
\end{definition}

\begin{lemma}
\label{lem:hoeffding}
Let $\left\{X_i\right\}_{i \geq 1}$  be i.i.d random variables defined on $\zp$, and $\E{X_1}<\ity$. For any $\delta \in (0,1)$, there exist $a,b >0$ such that 
\[\pb \left(\sum_{i=1}^n X_i \leq n(1-\delta)\E{X_1}\right) \leq a e^{-b\delta^2 n},
\]
for all $n\geq 1$.
\end{lemma}
\proof{Proof.}
Because $\pb\left(X_1 < 0\right)=0$ and $\E{X_1}<\ity$, by the dominated convergence theorem, 
\[\lim_{K\rightarrow \ity}\E{X_1 \cdot \mathbb{I}\left(X_1 \leq K\right)} = \E{X_1}.\]
Let $K(\beta) = \min \left\{K \in \zp: \E{X_1 \cdot \mathbb{I}\left(X_1 \leq K\right)}  \geq \beta \E{X_1}\right\}$ for some $\beta\in (0,1)$, and let $\left\{\tilde{X}^\beta_i\right\}_{i\geq 1}$ be a sequence of i.i.d random variables with distribution 
\[\tilde{X}^\beta_1 \stackrel{(d)}{=} X_1 \cdot \mathbb{I}\left(X_1 \leq K(\beta)\right).\]
Fix $\delta \in (0,1)$, and let $\beta = \frac{1-\delta}{1-\delta/2}$. By Hoeffding's inequality for bounded i.i.d random variables, there exist $a,b >0$ such that
\begin{eqnarray*}
\pb \left(\sum_{i=1}^n X_i \leq n(1-\delta)\E{X_1}\right) 
&\leq& 	\pb \left(\sum_{i=1}^n \tilde{X}^\beta_i \leq n(1-\delta)\E{X_1}\right) \\
&\stackrel{(a)}{\leq}& \pb \left(\sum_{i=1}^n \tilde{X}^\beta_i \leq n\frac{1}{\beta}(1-\delta)\E{\tilde{X}^\beta_1}\right)\\
&=& \pb \left(\sum_{i=1}^n \tilde{X}^\beta_i \leq n(1-\delta/2)\E{\tilde{X}^\beta_1}\right)\\
&\leq& a e^{-b \delta^2n},
\end{eqnarray*}
for all $n\geq 1$, where $(a)$ follows from the definition of $K(\beta)$ and $\tilde{X}^\beta_i$. 
\Halmos\endproof

\section{Other Proofs}

\subsection{Proof of Lemma \ref{lem:supadd}}
\label{app:lem:supadd}

\proof{Proof.} The scale-invariance property follows directly from the definition of $\phir{\cdot}$.
Fixing any $\bls\in\phir{\blw}$ and $\blts\in\phir{\bltw},$ we have that $0\leq\bls_{i}+\blts_{i}\leq\blw_{i}+\bltw_{i},$ and that
\[
\nor{\bls+\blts}=\nor{\bls}+\nor{\blts}\leq \norp{\blw}+\norp{\bltw}=\norp{\blw+\bltw}.
\]
This shows that $\bls+\blts\in\phir{\blw+\bltw},$ which proves the superaddativity. Finally, the convexity is a consequence of the first two properties: $\phir{a\blw+(1-a)\bltw} \supset \phir{a\blw}+ \phir{(1-a)\bltw} = a \phir{\blw}+(1-a) \phir{\bltw}.$
\Halmos\endproof

\subsection{Proof of Lemma \ref{lem:discount_Bell}}
\label{app:lem:discount_Bell}
\proof{Proof.} The lemma is a close analogue of the well-known Bellman's equation in a finite state space (cf.~\cite{Ross83}). We begin by converting the action set $\phis{\blw}$ into one that does not depend on the value of $\blw$, by constructing the following dynamic system. 
\begin{definition}
\textbf{SIM2} is a discrete-time dynamic system with\end{definition}
\begin{enumerate}
\item \textbf{States}: $\blw(t)\in\simp,\, t\in \zp$.
\item \textbf{Actions}: At time $t$, $\pi_{t}$ chooses a point $\bls\left(t\right)$
in $\simp$. 
\item \textbf{Transition}: $\blw(t+1)=L\left(p\left(\blw\left(t\right),\bls\left(t\right)\right)\right)$,
where $L\left(\bls\right) \bydef \proj{\calr\bls}$, and $p\left(\cdot,\cdot\right)$
is defined to be the projection of $\bls$ onto $\phis{\blw}$,%
\footnote{Note that the $\arg\min$ operation here always maps to a singleton
set, since $f\left(\bltw\right)\bydef\nor{\bltw-\bls}_{2}$ is a strictly
convex function in $\bltw$ for any $\bls,$ and the constraint set
$\phis{\blw}$ is convex and closed.%
} 
\[
p\left(\blw,\bls\right)\bydef\arg\min_{\bltw\in\phis{\blw}}\nor{\bltw-\bls}_{2},
\] 
where $\nor{\cdot}_{2}$ is the $L_{2}$ norm on $\R^K$ . 
\item \textbf{Reward per-stage}: $R_{2}\left(\blw(t),\bls(t)\right)=R\left(\blw(t), p\left(\blw\left(t\right),\bls\left(t\right)\right)\right)$.
\end{enumerate}
Compared to SIM, the action set at each stage in SIM2
is simply the entire $\simp$, while in the transition step, we first
project the chosen action, $\bls,$ to the feasible set $\phis{\blw}$,
before applying the reproduction matrix $\calr$. Given a current
state $\blw$ and a chosen action $\bls$ in SIM2, define
the transition map $T\left(\blw,\bls\right)\bydef L\left(p\left(\blw,\bls\right)\right)$. 
It is not difficult to verify from definitions that
\[
\max_{\bls\in\phis{\blw}}\left[R\left(\blw, \bls\right)+\gamma V^{\gamma}\left(\proj{\calr\bls}\right)\right]=\max_{\bls\in\simp}\left[R_2\left(\blw, \bls\right)+\gamma V^{\gamma}\left(T\left(\blw,\bls\right)\right)\right],
\]
for all $\blw\in\simp$. Hence, it suffices to show that 
\[
V^{\gamma}\left(\blw\right)=\min_{\bls\in\simp}\left[R_2\left(\blw, \bls\right)+\gamma V^{\gamma}\left(T\left(\blw,\bls\right)\right)\right],\quad\forall\blw\in\simp.
\]
Suppose that there exists an optimal stationary policy that achieves
the infimum in Eq.~\eqref{eq:vinf} for all $\blw\in\simp$, Lemma \ref{lem:discount_Bell} would follow directly by a simple inductive argument. To show the existence
of an optimal stationary policy, we invoke the following known result,
which is a special case for the more general statement given in \cite{Mai68}.
\begin{lemma} 
\label{lem:exsStat}For any discounted factor $\gamma\in\left(0,1\right),$
SIM2 admits an optimal stationary policy if the
following hold:
\begin{enumerate}
\item The action set, $\simp$, is a compact metric space.
\item The reward function $R_2\left(\cdot,\cdot\right)$ is bounded and
upper-semi continuous on $\simp\times\simp$. 
\item The transition map $T\left(\cdot,\cdot\right)$ is continuous on $\simp\times\simp$.%
\footnote{This condition is stated in a more general way in \cite{Mai68}. In
particular, letting $\delta_{a}\left(\cdot\right)$ be the Dirac measure
defined on $\simp$ with a unit mass on $a$, the condition requires
that for all $\blw,\bls\in\simp,$ if $\blw_{n}\rightarrow\blw$
and $\bls_{n}\rightarrow\bls$, then $\delta_{T\left(\blw_{n},\bls_{n}\right)}\left(\cdot\right)$
converges weakly to $\delta_{T\left(\blw,\bls\right)}\left(\cdot\right)$
as $n\rightarrow\infty$. One can verify that this is implied by the
continuity of $T\left(\cdot,\cdot\right)$.%
} 
\end{enumerate}

\end{lemma}
We now verify each of the three conditions in Lemma \ref{lem:exsStat}.
Condition $1$ and the boundedness of $R_2\left(\cdot,\cdot\right)$
in Condition $2$ follow directly from the definitions of $\simp$
and $R\left(\cdot,\cdot\right)$, respectively. 
Via some elementary analysis, one can show that $\phis{\blw}$ is convex-valued
and continuous in $\blw$ (i.e., both upper and lower semicontinuous), and that the function $p\left(\cdot,\cdot\right)$
is continuous; together, they imply that $R_2\left(\cdot,\cdot\right)$
is upper semicontinuous. Finally, the continuity of $T\left(\cdot,\cdot\right)$
follows from the continuities of $p\left(\cdot,\cdot\right)$ and
$L\left(\cdot\right)$. 
This completes the proof of Lemma \ref{lem:discount_Bell}. \Halmos\endproof

\subsection{Proof of Lemma \ref{lem:tech2}}
\label{app:tech2}
\proof{Proof.} Fix $\gamma\in(0,1)$. Let $A(\gamma,n) = \sum_{t=0}^{n-1} \gamma^t a_t$, and $B(\gamma,n) = \sum_{t=0}^{n-1} \gamma^t b_t$. Since $A(n,1)-B(n,1)\geq 0$ for all $n$ by assumption, it suffices to show that
\begin{equation}
	A(\gamma,n)-B(\gamma,n) \geq \gamma^{n}\left(A\left(1,n\right) - B\left(1,n\right)\right),
	\label{eq:seq1}
\end{equation}
for all $n\geq 1$. We do so by induction. Eq.~\eqref{eq:seq1} clearly holds for $n=1$. Suppose that it holds for some $n\geq 1$. We have
\begin{eqnarray*}
	&&A(\gamma,n+1)-B(\gamma,n+1) \\
	&=& \left(A(\gamma,n)-B(\gamma,n)\right) + \gamma^n\left( a_n - b_n\right)\\
	&\geq& \gamma^n \left(A(1,n)-B(1,n)\right) + \gamma^n\left( a_n - b_n\right)\\
	&\geq& \gamma^{n+1}\left(A(1,n+1)-B(1,n+1)\right),
\end{eqnarray*}
where the first inequality follows from the induction hypothesis, and the second inequality from the fact that $\gamma\in(0,1)$ and $A(1,n+1)-B(1,n+1)\geq0$. 
\Halmos\endproof

\subsection{Proof of Lemma \ref{lem:avgBell}}
\label{app:lem:avgBell}
\proof{Proof.} Let $\left\{ \left(\blw_{i},\bls_{i}\right)\right\} _{0\leq i\leq t}$ be a feasible sequence produced by some policy $\pi$. We have
\begin{eqnarray*}
g\left(\blw(t)\right) &=& g\left(\proj{\calr \bls\left(t-1\right)}\right)\\
& = & R\left(\blw(t-1),\bls(t-1)\right)+g\left(\proj{\calr\bls(t-1)}\right)-R\left(\blw(t-1),\bls(t-1)\right)\\
 & \leq & \left(\max_{\bls\in\phis{\blw(t-1)}}R\left(\blw(t-1),\bls\right)+g\left(\proj{\calr\bls}\right)\right)-R\left(\blw(t-1),\bls(t-1)\right)\\
 & = & \alpha+g\left(\blw(t-1)\right)-R\left(\blw(t-1),\bls(t-1)\right),
\end{eqnarray*}
where the last equality follows from the definition of $g$ in Eq.~\eqref{eq:avgBell},
with equality holding for $\pi=\pi^{*}$. This implies that, for all $n>0$,
\[	\sum_{t=1}^{n}\left[g\left(\blw(t)\right)-\left(\alpha+g\left(\blw(t)\right)-R\left(\blw(t-1),\bls(t-1)\right)\right)\right] \leq 0,
\]
or, equivalently, that
\[
\alpha\geq\frac{g\left(\blw(t)\right)}{n}-\frac{g\left(\blw(t)\right)}{n}+\frac{\sum_{t=1}^{n}R\left(\bls(t-1)\right)}{n},
\]
with equality holding if $\pi=\pi^*$. The claim follows by taking the limit of $n\rightarrow \infty$, and noting that both $g(\cdot)$ and $R\left(\cdot,\cdot\right)$ are bounded over $\simp$. 
\Halmos\endproof

\subsection{Proof of Lemma \ref{lem:boundary}}
\label{app:lem:boundary}
\proof{Proof.} 
Let $\bls \in \phis{\blw}\backslash \partial\left(\phis{\blw}\right)$, and $\proj{\calr \bls} = \blw$. We will show that $\blw$ must not be an optimal fixed point via a perturbation argument. Let $\bls(d) = \bls + (d , -d)^\top$. Since $\bls$ is in the interior of $\phis{\blw}$ and the scaling map $\proj{\cdot}$ is continuous, there exists $d>0$ such that 
\[\bls(d) \in \phis{\blw}, \mbox{ and } \bls\in \proj{\calr \bls(d)}.\]
Therefore, the state-action sequence constructed by alternating between $\left(\blw,\bls(d)\right)$ and $\left(\proj{\calr\bls(d)},\bls\right)$ is feasible, and its average reward is given by $\frac{1}{2}\beta\left(\ln\nor{\calr\bls}+ \ln\nor{\calr\bls(d)}\right)$ which is strictly greater than that of the original sequence, $\beta\ln \nor{\calr\bls} $, by noting our assumption that $\nor{\calr_{(\cdot,1)}}>\nor{\calr_{(\cdot,2)}}$. This proves our claim that the corresponding subpopulation for an optimal fixed-point $\blw$ must lie on the boundary of $\phis{\blw}$. 
\Halmos\endproof

\section{Discrete-time Cancer Growth Model}
\label{app:cancer}

In this section, we derive a discrete-time multi-type branching process from a continuous-time active-quiescent cell-cycle kinetics model for cancer growth. Denote by $\blw(t)=\left(\blw_1(t),\blw_t(2)\right)^\top$ the cell population in round $t$ of the process, where $\blw_1(t)$ and $\blw_2(t)$ are the numbers of active and quiescent cells, respectively, and two consecutive rounds are 3 weeks apart. Assuming \emph{no treatment} is applied, we would like to show that 
\begin{equation}
	\blw(t+1)=\calr\blw(t),
	\label{eq:w_discrete}
\end{equation}
where the value of $\calr$ is given in Eq.~\eqref{eq:cancer_matrix}. 

We use a continuous-time dynamic model proposed in \cite{DH06}. Let $s \in \R_+$, and as before, denote by $\blw(s)$ of the cell population at time $s$. The following system of differential equations governs the evolution of $\left\{\blw(s):s\in \R_+\right\}$:
\begin{equation}
	\dot{\blw}(s) = A\blw(s), \quad t\in \R_+, 
	\label{eq:diff_w}
\end{equation}
where 
	\[
	A = \left( \begin{array}{cc} -\mu & \, \gamma \\ 2\mu  & \, -\gamma \end{array} \right).
\]
The biological process underlying Eq.~\eqref{eq:diff_w} is as follows: active cells undergo mitosis at rate $\mu$, and when this occurs, two daughter cells are produced and directly enter the quiescent state; quiescent cells do not reproduce, but re-enter the active compartment at rate $\gamma$. The parameter values used in \cite{DH06} are: 
	\[
	\mu=0.0655, \mbox{ and } \gamma = 0.0476,
\]
and a unit of $s$ corresponds to one day. 

Assume that $A$ is an invertible $K\times K$ matrix. It is well-known that Eq.~\eqref{eq:diff_w} admits a unique solution of the following form. Denote by $V = \left(\blv_1,\ldots,\blv_K\right)$ the matrix where $\blv_i$ is the $i$th right eigenvector of $A$, and by $\lambda_i$ the  eigenvalue associated with $\blv_i$. Let $D(s)$ be the diagonal matrix where $D_{i,i}(s) = e^{\lambda_i t}$ for all $i\in \left\{1,\ldots, K\right\}$. Then the solution to Eq~\eqref{eq:diff_w} is given by 
\begin{equation}
	\blw(s) = V D(s) V^{-1}\blw(0), \quad s \in \R_+.
	\label{eq:diff_sol}
\end{equation}
By the Eq.~\eqref{eq:diff_sol}, it is easy to check that the process  $\left\{\blw(s):s\in \R_+\right\}$ sampled at points $s_t=21 t$, $t\in \zp$, is a discrete-time branching process, and the reproduction matrix $\calr$ (Eq.~\eqref{eq:w_discrete}) is obtained by setting 
	\[
	\calr =  V D(21) V^{-1}.
\]
\end{APPENDIX}


\begin{thebibliography}{References}
\footnotesize

\bibitem{Mai68}A.\ Maitra, 
``Discounted dynamic programming on compact metric spaces,''
\textit{Sankhya, Ser. A}, 27:241--248, 1968.

\bibitem{Ross83}S.\ Ross, ``Introduction to Stochastic Dynamic Programming,''
\textit{Academic Press}, 1983.

\bibitem{AG05} L.\ A.\ Adamic and N.\ Glance,
``The political blogosphere and the 2004 U.S. election: Divided they blog,''
\emph{Proceedings of KDD Workshop on Link Analysis and Group Detection}, LinkKDD, 2005. 

\bibitem{SPK96}
A.\ Swierniak, A.\ Polanski, and M.\ Kimmel. 
``Optimal control problems arising in cell-cycle-specific cancer chemotherapy.''
\emph{Cell Proliferation,} 29:117-139, 1996.

\bibitem{DH06}
A.\ Dawson and T.\ Hillen,
``Derivation of the tumour control probability (TCP) from a cell cycle model,''
\emph{Computational and Mathematical Methods in Medicine}, 7:121-141, 2006.

\bibitem{He84}
G.\ H.\ Heppner, ``Tumor heterogeneity,'' 
\emph{Cancer Research}, 44:2259-2265, 1984.

\bibitem{PC07}
K.\ Polyak and L.\ Campbell, ``Breast tumor heterogeneity: Cancer stem cells or clonal evolution?''
\emph{Cell Cycle}, 6(19):2332-2338, 2007.

\bibitem{DF11}
R.\ Durrett, J.\ Foo, K.\ Ledeer, J.\ Mayberry, and F.\ Michor, ``Intratumor heterogeneity in evolutionary models of tumor progression,''
\emph{Genetics}, 188:461-477, 2011. 

\bibitem{DDHM12}
K.\ Danesh, R.\ Durrett, L.\ Havrilesky, and E.\ Myers,
``A branching process model of ovarian cancer,''
\emph{Journal of Theoretical Biology}, 314:10 - 15, 2012. 

\bibitem{Har89}
T.\ E.\ Harris, 
``The theory of branching processes,''
New York: Dover Publishers, 1989. 

\bibitem{SZ74}
B.\ A.\ Sevastyanov and A.\ M.\ Zubkov,
``Controlled branching processes,''
{\em Theory Probab. Appl.}, 19:14--24, 1974.

\bibitem{Pli77}
S.\ R.\ Pliska,
``Optimization of multitype branching processes,''
\emph{Management Sci.}, 23(2):117--24, 1976. 

\bibitem{RW82}
U.\ G.\ Rothblum and P.\ Whittle,
``Growth optimality for branching Markov decision chains,''
\emph{Math. of Operations Research}, 7:582--601, 1982.

\bibitem{GP02}
M.\ Gonzalez, M.\ Molina, and I.\ del Puerto 
``On the class of controlled branching process with random control functions,''
\emph{J. Appl. Probab.}, 40:995--1006, 2002.
 

\bibitem{LBEW10}
R.\ V.\ der Lans, G.\ van Bruggen, J.\ Eliashberg, and B.\ Wierenga ``Viral branching model for predicting the spread of electronic word-of-mouth,'' \emph{Marketing Science}, 29(2):348--365, 2010.

\bibitem{IM09}
J.\ L.\ Iribarren and E.\ Moro,
``Impact of human activity patterns on the dynamics of information diffusion,''
\emph{Physical Rev. Lett.}, 103(3), 2009.

\bibitem{Vaz2006}
A.\ Vazquez,
``Spreading dynamics on heterogeneous populations: multi-type network approach,''
\emph{Physical Review Letters}, 74, 2006.

\bibitem{RD02}
M.\ Richardson and P.\ Domingos,
``Mining knowledge-sharing sites for viral marketing,''
\emph{8th Intl. Conf. on Knowledge Discovery and Data Mining}, 2002.

\bibitem{KKT03}
D.\ Kempe, J.\ Kleinberg, and E.\ Tardos,
``Maximizing the spread of influence through a social network,'' 
\emph{Proc. KDD}, 2003.

\bibitem{FKG03}
C.\ Farrington, M.\ Kanaan and N.\ Gay,
``Branching process models for surveillance of infectious diseases controlled by mass vaccination,''
\emph{Biostatistics}, 4:279--295, 2003.

\bibitem{Res93}
J.\ A.\ C.\ Resing, 
``Polling systems and multitype branching processes,'' 
\emph{Queuing Syst.}, 33(4):. 409-426, 1993.

\bibitem{Ok07}
E.\ Ok, 
\emph{Real Analysis With Economic Applications,}
Princeton Press, 2007. 

\bibitem{Kak41}
S.\ Kakutani,
``A generalization of Brouwer's fixed point theorem,''
\emph{Duke Mathematical Journal}, 8(3): 457-459, 1941.


\end{thebibliography}
\end{document}